\documentclass[11pt,a4paper,twoside]{amsart}
\usepackage{amssymb,latexsym,amsmath,graphics}
\usepackage{graphicx}
\usepackage{caption}
\usepackage{subcaption}
\usepackage{enumerate}  
\usepackage{enumitem}
\usepackage{color}

\newlist{abbrv}{itemize}{1}
\setlist[abbrv,1]{label=,labelwidth=1in,align=parleft,itemsep=0.1\baselineskip,leftmargin=!}

\numberwithin{equation}{section}

\newcommand{\ntau}{\tau }

\newcommand{\mmattitext}{}

\newcommand{\mntext}{}
\newcommand{\ntext}{}
\newcommand{\newtext}{}

\newtheorem {theorem}{Theorem}
\newtheorem {lemma}{Lemma}

\newtheorem {definition}{Definition}

\theoremstyle{definition}
\newtheorem{remark}{Remark}
\newtheorem {example}{Example} 

\newcommand{\argmin}{\mathop{\mathrm{arg\,min}}}

\newcommand{\finm}{\mathbf{m}}
\newcommand{\finA}{\mathbf{A}}

\newcommand{\T}{{\mathbb T}}
\newcommand{\R}{{\mathbb R}}
\newcommand{\E}{{\mathbb E}}
\newcommand{\N}{{\mathbb N}}
\newcommand{\Z}{{\mathbb Z}}

\newcommand{\F}{{\mathcal F}}

\newcommand{\ra}{\rightarrow}
\newcommand{\expec}{{\mathbb E \,}}

\newcommand{\beqno}{\begin{eqnarray*}}
\newcommand{\eeqno}{\end{eqnarray*}}

\newcommand{\beqla}[1] {\begin {eqnarray}\label{#1}}

\def \beq {\begin {eqnarray}}
\def \eeq {\end {eqnarray}}
\def \ba {\begin {eqnarray*}}
\def \ea {\end  {eqnarray*}}

\def \b {\beta}
\def \a {\alpha}

\def \bra{\langle}
\def \cet{\rangle}
\def \prob{{\mathbb P}}

\def \Noise {\mathcal{E}}

\title[\resizebox{4.5in}{!}{Posterior consistency and convergence rates for Bayesian inversion}]{Posterior consistency and convergence rates for Bayesian inversion with hypoelliptic operators}

\author{Hanne Kekkonen, Matti Lassas and Samuli Siltanen}
\thanks{Department of Mathematics and Statistics
P.O. Box 68,
00014 University of Helsinki
FINLAND, e-mails: hanne.kekkonen@helsinki.fi, matti.lassas@helsinki.fi, samuli.siltanen@helsinki.fi}
\address {Department of Mathematics and Statistics
P.O. Box 68,
00014 University of Helsinki
FINLAND}

\begin{document}

\maketitle

\begin{abstract}
Bayesian approach to inverse problems is studied in the case where the forward map is a linear hypoelliptic pseudodifferential operator and measurement error is additive white Gaussian noise. The measurement model for an unknown Gaussian random variable $U(x,\omega)$ is 
\begin{eqnarray*}
M_\delta(y,\omega) = A(U(x,\omega) )+ \delta\hspace{.2mm}\mathcal{E}(y,\omega),
\end{eqnarray*}
where $A$ is a finitely many orders smoothing linear hypoelliptic operator and $\delta>0$ is the noise magnitude. The covariance operator $C_U$ of $U$ is smoothing of order $2r$, self-adjoint, injective and elliptic pseudodifferential operator.

If $\mathcal{E}$ was taking values in $L^2$ then in Gaussian case solving the conditional mean (and maximum a posteriori) estimate is linked to solving the minimisation problem 
\begin{eqnarray*}
T_\delta(m_\delta) = \argmin_{u\in H^r}
  \big\{\|A u-m_\delta\|_{L^2}^2+ \delta^2\|C_U^{-1/2}u\|_{L^2}^2 \big\}. 
\end{eqnarray*}
However, Gaussian white noise does not take values in $L^2$ but in $H^{-s}$ where $s>0$ is big enough. A  modification of the above approach to solve the inverse problem is presented, covering the case of white Gaussian measurement noise. Furthermore, the convergence of conditional mean estimate to the correct solution as $\delta\rightarrow 0$ is proven in appropriate function spaces  using microlocal analysis. Also the frequentist posterior contractions rates are studied. 
\end{abstract}  

\noindent\keywords{{\bf Keywords:} Posterior consistency, convergence rate, Bayesian inverse problem, white noise, pseudodifferential operator}
 
\tableofcontents

\section{Introduction}

\noindent
Practical inverse problems arise from the need to extract information from indirect data. For example, consider a 
device designed for measuring point values of a physical quantity $u(x)$.  Technological imperfections cause the values of $u$ at nearby points to merge together in the measurement. Mathematically this corresponds to convolution $\Phi\ast u$ by a point spread function $\Phi$. The inverse problem is to recover the function $u$ approximately from a finite number of point values of $\Phi\ast u$ corrupted by random white noise.  

Computational inversion requires a finite representation of the quantity $u(x)$. In this paper we promote the view that it is a good idea to design a continuous model for $u$ which then can be discretized with a desired number of degrees of freedom. This paves the way for the analysis of convergence as the discretization becomes finer. Such convergence enables switching between different discretizations consistently; this is crucial for multigrid methods and for certain parameter-choice strategies.

In ill-posed inverse problems the measurement data alone is not sufficient for noise-robust recovery of the quantity of interest. For instance, Fourier transforming $\Phi\ast u$ gives $\widehat{\Phi}\widehat{u}$, so frequency-domain information is lost in areas where $\widehat{\Phi}$ is very close to zero. Therefore, successful computational inversion requires some {\em a priori} information in addition to the measurement data. 

Practical inversion is all about combining measurement data and {\em a priori} information in a noise-robust way. The classical approach to do this is {\em regularization} that assumes that the noise is deterministic and small in norm. Regularization involves defining a family of continuous maps, parametrized by the norm $\|{\rm noise}\|_{L^2}$, from the data space to the space of unknown quantities. This must be done so that as $\|{\rm noise}\|_{L^2}\to 0$, the reconstruction approaches the true solution along a stable path. 
This methodology was originated by Tikhonov \cite{Tikhonov1943,Tikhonov1963}. Both continuous and discrete cases have been studied in depth in \cite{Engl1996a,Groetsch1984,Kirsch1996,Phillips1962,Morozov1984, Tikhonov1995, schuster2012}. 

There is a serious drawback in the above noise model in the continuous limit. Namely, continuous white noise on $\R^d$ is not square integrable. We discuss this in detail below in Section \ref{sec:whitenoise}. The goal of this paper is to use  {\em Bayesian inversion} to construct a consistent continuous-discrete framework covering the case of white noise.  

Bayesian inversion is a flexible framework for combining measurement data and  {\em a priori} information in the form of a {\em posterior distribution} \cite{Cavalier2008, Fitzpatrick1991, Kaipio2004a, Lassas2009, Stuart2010}. 
Computational exploration of the finite-dimensional posterior distribution yields useful estimates of the quantity of interest and enables uncertainty quantification. Furthermore, analytic results about the continuous model can then be restricted to a given resolution in a discretization-invariant way.

Our approach to Bayesian inversion follows a general strategy of computational mathematics: we consider a continuous model which can be discretized for any practical setting. 

We study the following continuous model for indirect measurements:
\begin{equation}\label{infinitemodel1}
  M_\delta = AU + \Noise\delta,
\end{equation} 
where the  random variables $M_\delta$ (data) and $U$ (quantity of interest) take values in the Sobolev spaces $H^{-s}(N)$ and $H^{\ntau}(N)$, respectively.  Here $N$ is a $d$-dimensional compact  manifold e.g. a torus corresponding to a $d$-dimensional cube with opposite sides glued together. The real parameter $\ntau$ is related to our {\em a priori} information about the smoothness of the unknown quantity of interest.

The measurement operator $A$ in our model (\ref{infinitemodel1}) is quite general: we assume it to be a finitely smoothing, injective hypoelliptic pseudodifferential operator ($\Psi$DO). See section \ref{sec:generalproof} for precise definition. This class includes convolution operators with finitely smooth kernel $\Phi$.  One example of an operator that is hypoelliptic but not elliptic is the heat operator. For more examples of hypoelliptic operators see Appendix \ref{ap:hypo}. The measurement noise $\Noise$ is assumed to be normalised white Gaussian noise with mean zero and unit variance, and $\delta>0$ models the noise amplitude. 

We model practical measurement data by 
\begin{equation}\label{eq: measmodel_k k}
  M_k = P_k(A U)+P_k(\Noise)\delta.  
\end{equation}
Here $P_k$ is a linear operator related to measurement device; we assume that $P_k$ is an orthogonal projection with $k$-dimensional range. We discretize the unknown $U$ using some computationally feasible approximation of the form $U_n=S_n U$. Now we can study an inverse problem
\begin{equation}\label{discreteinverseproblem}
\begin{split}
\mbox {\em given a realisation of $M_k$, estimate $U_n$}. 
\end{split}
\end{equation}

We are interested to know what happens to the approximated solutions of (\ref{infinitemodel1}) when $\delta\to 0$. The analysis of small noise limit, also known as the theory of posterior consistency, has attracted a lot of interest in the last decade. Posterior convergence rates were first studied in \cite{Ghosal2000,Shen2001} and further in papers \cite{Agapiou2013, Agapiou2014, Dashti2013, Hoffmann2015, Huang2004, Knapik2014, Knapik2012, Knapik2013, Ray2013, Shen2001, szabo2015a, Vollmer2013}.  
However, much remains to be done.
Developing a comprehensive theory is important since posterior consistency justifies the use of the Bayesian approach the same way as convergence results justify the use of regularisation techniques.

In the above mentioned papers the problem is studied from the frequentist point of view, that is, the data is thought to be generated by a fixed 'true' solution $u^\dagger$ instead of random draw $U(\omega)$ from the prior distribution. This means that all the randomness in $M$ comes from the randomness of the noise $\Noise$. The interest is then on the contraction of the posterior distribution around the 'true' solution $u^\dagger$ as the noise goes to zero, see Subsection \ref{sec: frequentist} and Theorem  \ref{Thm:contraction}. The main emphasis of this paper, however, is in the purely Bayesian approach where it is assumed that also $U$ is random. Since $U$ and $\Noise$ are assumed independent we can write
\begin{equation}
  M_\delta(\omega) = AU(\omega_1) + \Noise(\omega_2)\delta,
\end{equation}   
where $\omega=(\omega_1,\omega_2)\in \Omega_1\times\Omega_2$. In Bayesian case the posterior distribution is a function of $(\omega_1,\omega_2)$. Also the probability measure $d\mathbb{P}$ can be written in the following form 
\begin{align*}
d\mathbb{P}=d\mathbb{P}_1(\omega_1)d\mathbb{P}_2(\omega_2). 
\end{align*}
We will denote the expected value  over the joint distribution of $U$ and $\Noise$ by $\E$.  The expected value over the noise is defined by
\begin{align}\label{frequentis_mean}
\E_U(F(\omega_1,\omega_2))=\int_{\Omega_2}F(\omega_1,\omega_2)d\mathbb{P}_2(\omega_2).
\end{align}

Our paper provides a conceptual advantage over much of the existing literature. In many earlier studies $A$ and $C_U$ are perturbations of negative powers of operator $(I-\Delta)$. Our assumption, formulated in terms of hypoelliptic operators, means roughly speaking that the measurement operator $A$ and the covariance operator $C_U$ do not need to have a common basis in their singular value decomposition.  

The rest of this paper is organised as follows. In section \ref{sec:mainresult} we will introduce the Bayesian setting we are using and present our main result about convergence rates. In section \ref{sec:Gaussianexample} we take a closer look to the generalised Gaussian random variables in Sobolev spaces and introduce so called white noise paradox. We will also show that the distribution $U$ takes values in $H^{\ntau}$, where $\ntau$ is related to the smoothness of the solution and depends on the dimension of the space and the covariance of the prior. In section \ref{sec:generalproof} we will introduce hypoelliptic operators and prove Theorems \ref{theorem:main} and \ref{thm:frequentist1}. In section \ref{sec:crediblesets} we characterise credible sets and frequentist confidence regions and present and prove two theorems about the contraction of them. In Appendix \ref{ap:hypo} we will give examples of some hypoelliptic operators and in Appendix \ref{ap: computational} we give a computational example.

\subsection*{Notations}

\begin{abbrv}
\item[$S^m$]		Class of pseudodifferential symbols of order $m$. See Definition \ref{def: pseudo_symbol}.
\item[$\Psi^m$]		Space of pseudodifferential operators ($\Psi$DO) of order $m$. See Definition \ref{def: pseudo_operator}.
\item[$H\Psi^{m,m_0}$]		Space of hypoelliptic $\Psi$DOs of type $(m,m_0)$. See Definition \ref{def: hypoelliptic}.
\item[$\Psi^{m}_p$]  Space of  $\Psi$DOs of order $m$ depending on spectral variable with order $p$. See Definition \ref{def: spectral}.
\item[$H\Psi^{m,m_0}_p$] Space of hypoelliptic $\Psi$DOs of type $(m,m_0)$ depending on spectral variable with order $p$. See Definition \ref{def: spectral_hypo}.
\item[$Tr_{H^q\to H^q }C$]		The trace of the operator $C:H^q\to H^q$. See (\ref{tracecondition}).

\end{abbrv}

\section{Convergence results}\label{sec:mainresult}

Let us return to our indirect measurement problem 
\begin{equation*}
  M_\delta = AU + \Noise\delta,
\end{equation*}  
where we model $U=U(x,\omega),\ M_\delta=M_\delta(y,\omega)$ and $\Noise=\Noise(y,\omega)$ as random functions.
Here $\omega\in\Omega$ is an element of a complete probability space  $(\Omega,\Sigma,\prob)$ and $x$ and $y$ denote the variables in domains of Euclidean spaces. The reason why we model also $U$ as a random variable is that even though the unknown quantity is assumed to be deterministic we have only incomplete data of it. All information available about $U$ before performing the measurements is included in a prior distribution that is independent of the measurement.

The Bayesian inversion theory is based on the Bayes formula.
To solve the inverse problem (\ref{discreteinverseproblem}) we have to express available a priori information of $U_n$ in the form of a prior distribution $\pi_{pr}$ in an $n$-dimensional subspace. Let $M_k$ and $\Noise_k$ be random vectors taking values in $\R^k$, and denote their distributions by $\pi_{M_k}$ and $\pi_{\Noise_k}$, respectively. The solution of the inverse problem after performing the measurements is the posterior distribution of the unknown random variable.  
Given a realisation of the discrete measurement the posterior density for $U_n$ taking values in the $n$-dimensional subspace is given by the Bayes formula 
\begin{eqnarray}\label{posterior2}   
\begin{split}  
\pi(\textbf{u}\,|\,\textbf{m}_\delta) &=\frac{\pi_{pr}(\textbf{u})\pi_{\Noise_k}(\textbf{m}_\delta\ |\ \textbf{u})}{\pi_{M_k}(\textbf{m}_\delta)}\\
& =  c\pi_{pr}(\textbf{u})\exp\Big(-\frac{1}{2\delta^2}\|\textbf{m}_\delta - \textbf{A} \textbf{u}\|_2^2\Big),\quad \textbf{u}\in\R^n,\ \textbf{m}_\delta\in\R^k
\end{split}
\end{eqnarray}
where $\textbf{A}=P_kAS_n$ is a $k\times n$ matrix approximation to the operator $A$.

An approximated solution for the inverse problem is often given as a point   estimate for (\ref{posterior2}).
Let us assume that also $U_n$ has Gaussian distribution. The maximum a posteriori (MAP) estimate $T^{MAP}_\delta: \R^k\to\R^n$ is defined by
\begin{equation}\label{mapkaava}
T^{MAP}_\delta(M_k(\omega)):=\arg\max_{\textbf{u}\in \R^n}\pi(\textbf{u}\,|\,M_k(\omega)).
\end{equation}
Note that the MAP estimate depends on $\omega$ through the realisation of the noise $\Noise_k(\omega)$ and unknown $U_n(\omega)$.
When $U_n$ is Gaussian distributed the MAP estimate coincides almost surely with the conditional mean estimate (CM)
\begin{equation}\label{meankaava2}
  T^{CM}_\delta(M_k(\omega)) = \E(U_n | \mathcal{M}_k)(\omega) \quad \text{a.s.}
\end{equation}
where $\mathcal{M}_k$ is the $\sigma$-algebra generated by $M_k$.

Since in our case CM=MAP a.s. we will consider below the MAP estimate. Let us denote the covariance matrix of $U_n$ by $\textbf{C}_{U_n}$. Solving the maximisation problem (\ref{mapkaava}) with a fixed realisation of noise and unknown corresponds to solving the minimisation problem
\begin{equation}\label{discrTikh1}
T_\delta(\textbf{m}_\delta) = \mbox{arg}\min_{\textbf{u}\in \R^n}
  \big\{\frac 1{2\delta^2} \|\textbf{Au}-\textbf{m}_\delta\|_2^2+ \frac{1}{2}\|\textbf{C}_{U_n}^{-1/2}\textbf{u}\|_2^2 \big\}.
\end{equation}

Constructing $S_n$ and $\pi_{pr}$ is the core difficulty in Bayesian inversion. In {\mntext many inverse problems there is} no natural discretisation for the continuum quantity $U$, so $n$ can be freely chosen. Consequently, $S_n$ and $\pi_{pr}$ should in principle  be described for all $n>0$. This raises the following questions: do the chosen $S_n$ and $\pi_{pr}$ represent the same 
{\em a priori} knowledge consistently at all resolutions $n>0$? Does the estimate $T_\delta(\textbf{m}_\delta)$ converge as $n\ra\infty$? See e.g. \cite{kolehmainen2012,lassas2004}
Also, the number of data points may change, for example due to an updated measurement device.
The aim of this paper is to build a rigorous theory that allows us to connect discrete models to their
infinite-dimensional limit models in a consistent way.

We achieve consistent representation of {\em a priori} knowledge by constructing the prior
distribution for $U$ in the infinite-dimensional space $X$. Then the random variable $U_n=S_nU$
takes values in the finite-dimensional subspace $X_n\subset X$ and represents approximately the
same {\em a priori} knowledge {\ntext as} $U$. The same way we construct distributions for $M$ and $\Noise$ in the infinite-dimensional space $Y$ in which case the random variables $M_k$ and $\Noise_k$ take values in the finite-dimensional subspace $Y_k\subset Y$.

The finite-dimensional  problem (\ref{discrTikh1}) $\Gamma$-converges as $n,k\to \infty$, under certain assumptions (including that $m$ should be an $L^2$-function), to the following infinite-dimensional minimisation problem in a Sobolev space $H^r$:
\begin{equation}\label{contTikh0}  
 \argmin_{u\in H^r}\big\{\frac{1}{2\delta^2}\|m_\delta-A u\|_{L^2(N)}^2 + \frac{1}{2}\|C_U^{-1/2}u\|_{L^2(N)}^2 \big\}.
\end{equation}
Above $C_U^{-1/2}\in\Psi^r$, that is, $C_U^{-1/2}$ is $-r$ orders smoothing pseudodifferential operator.
See \cite{Kekkonen2014} for a proof.
If we are thinking the above as a MAP estimate to a Bayesian problem we have to assume that $U$ has formally the following distribution 
\begin{equation*}
\pi_{pr}(u) \underset{\mbox{\tiny \em formally}}= c\exp\bigg(-\frac{1}{2}\|C_U^{-1/2}u\|^2_{L^2(N)}\bigg).
\end{equation*}
Formula (\ref{contTikh0}) only makes sense if the noise is square integrable. Even though 
\ba
\|\varepsilon_k\|_{\ell^2}<\infty  
\ea
with any $k\in\N$ the limit, when $k\to \infty$, is infinity.
We will return to this so called `white noise paradox' in section \ref{sec:whitenoise}.

\subsection{Main result}  
Let us now modify formula (\ref{contTikh0}) to arrive at something useful for white Gaussian noise. 
When $\varepsilon\in L^2$ we can write 
\begin{equation}
\|m_\delta-A u\|_{L^2(N)}^2=\|A u\|_{L^2(N)}^2- 2( m_\delta,A u )_{L^2(N)} +\|m_\delta\|_{L^2(N)}^2.
\end{equation} 
Now omitting the infinite `constant term' $\|m_\delta\|_{L^2(N)}^2$  we get a new minimisation problem which is well-defined also when $m$ is not an $L^{2}$ function
\begin{equation}\label{contTikh1}
  T_\delta(m_\delta):=\argmin_{u\in H^r(N)} 
  \big\{\|A u\|_{L^2(N)}^2 - 2\langle m_\delta,A u \rangle +
 \delta^2\|C_U^{-1/2}u\|_{L}^2 \big\},
\end{equation}  
where $\langle m_\delta, A u \rangle$ is interpreted as a suitable duality pairing instead of $L^2(N)$ inner product. When {\newtext $A\in\Psi^{-t}$, $t\geq -\ntau+s$, we can define $\langle m_\delta, Au \rangle = \langle m_\delta, Au \rangle_{H^{-s}(N)\times H^{s}(N)}$. Note that the forward operator $A$, the prior distribution and the noise depend on on each other only through assumption $t\geq -\ntau+s$. 

It is well-known that the solution of the finite-dimensional problem (\ref{discrTikh1}) can be calculated using the following  formula:
\begin{equation}\label{discrTikh2}
  T_\delta(\textbf{m}_\delta)=(\finA^T\!\finA+\delta^2 \textbf{C}_{U_n}^{-1})^{-1}\finA^T \finm_\delta.
\end{equation}    
We can write the approximated solution $u_\delta:= T_\delta(m_\delta)$ of the continuous problem (\ref{contTikh1}) by
\begin{equation}\label{contTikh2}
T_\delta(m_\delta) =( A^*A+ \delta^2C_U^{-1})^{-1}A^*m_\delta.
\end{equation}

Before the main result of the paper we will study a simple example to give a reader an insight to Bayesian settings. 

\begin{example}
Let $N$ be a $1$-dimensional torus $\T^1$. We are interested of the inverse problem 
\ba 
M_\delta=AU+\delta\Noise
\ea
where we assume that $\Noise\sim N(0,I)$ and $U\sim N(0,I)$, that is both the noise and the unknown are assumed to be normalised white Gaussian noise, see section \ref{sec:Gaussianexample} for rigorous definition. The white noise takes values in $H^{\ntau}$ with some $\ntau<-1/2$. On the other hand white noise has formally the following distribution 
\begin{equation*}
\pi_{pr}(u) \underset{\mbox{\tiny \em formally}}= c\exp\bigg(-\frac{1}{2}\|u\|^2_{L^2(\T^1)}\bigg).
\end{equation*}
Hence we want to solve  
\begin{equation*}
 \argmin_{u\in L^2}\big\{\|A u\|_{L^2}^2 -2\langle m_\delta,Au\rangle+ \delta^2\|u\|_{L^2(\T^1)}^2 \big\}.
\end{equation*}
Note that we are looking for a solution in $L^2(\T^1)$ even though the realisations of $U$ are in $L^2(\T^1)$ with probability zero.
In general if we are interested in finding a solution in $H^r(N)$ then we can show that the prior should take values in $H^{\ntau}(N)$ where $\ntau=r-s$, see section \ref{sec:Gaussianexample}.
\end{example}

\subsubsection{Convergence results in Bayesian setting}
We will now formulate the main theorem of this paper about the convergence of the continuous solution (\ref{contTikh2}).  The precise definitions are discussed in more detail in section \ref{sec:Gaussianexample} and Theorem \ref{theorem:main} is proved in section \ref{sec:generalproof}.

\begin{theorem}\label{theorem:main}     
Let  $r,s\in [0,\infty)$ and $N$ be a $d$-dimensional closed manifold. Let $U(x,\omega)$ be a generalised Gaussian random function taking values in $H^{\ntau}(N)$, $\ntau=r-s$, with zero mean and covariance operator $C_U$. Assume that the operator $C_U$ is a self-adjoint, injective and elliptic pseudodifferential operator ($\Psi$DO) of order $-2r$. Let $\Noise(y,\omega)$ be white Gaussian noise on $N$. Consider the measurement 
 \begin{align*}
M_\delta(y,\omega) = A(U(\cdot,\omega))+\delta\Noise(y,\omega),\quad\quad \omega\in\Omega,
\end{align*}
where $A\in H\Psi^{-t,-t_0}$, $t>\max\{0, -\ntau+s\}$ and $t\leq t_0< 2t+r$, is a hypoelliptic pseudodifferential operator on the manifold $N$ and $A:L^2(N)\to L^2(N)$ is injective. Above $\delta\in \R_+$ is the noise level and $\Noise$ takes values in $H^{-s}(N)$ with some $s>d/2$.

Take $\zeta < \ntau-3(t_0-t)$. Then we have the following convergence
\begin{align}\label{zeta convergence}
\E\|U_\delta(\omega)-U(\omega)\|_{H^{\zeta}}\to 0,\quad \text{as}\,\,\, \delta\to 0.
\end{align}
The expectation is taken with respect to the joint distribution of $(U,\Noise)$.
We have the following estimates for the speed of convergence:
\begin{itemize}   
\item[(i)] If $\zeta\leq t-s-2t_0$ then there is such $C>0$ independent of $\delta$ that 
\begin{align}\label{converge_1}
\E\|U_\delta(\omega)-U(\omega)\|_{H^{\zeta}}\leq C\delta^{\frac{2t-t_0+r}{t_0+r}}.
\end{align}  
\item[(ii)] If $t-s-2t_0\leq \zeta< \ntau-3(t_0-t)$ then there is such $C>0$ independent of $\delta$ that 
\begin{align}\label{converge_2}
\E\|U_\delta(\omega)-U(\omega)\|_{H^{\zeta}}\leq C\delta^{-\frac{\zeta-\ntau+3(t_0-t)}{t_0+r}}.
\end{align} 
\end{itemize}   
\end{theorem}

The different convergence speeds (i) and (ii) show the trade-off between the smoothness of the space and the speed of convergence. In case (i) we get better convergence rates but in case (ii) we can use a stronger norm. We see also that the smoother the forward operator $A$ is the worse convergence rates we get.

We note that instead of the estimates (\ref{converge_1}) or (\ref{converge_2}), we could alternatively take the expected value only with respect to the noise in which case the constant $C$ would depend on the realisation of $U(\omega)$. That is, proof of Theorem 1 also shows that, we have almost surely  
\begin{align}\label{M1}
\limsup_{\delta\to 0} \frac{\E_U\|U_\delta(\omega)-U(\omega)\|_{H^{\zeta}}}{\delta^\nu}<\infty
\end{align}
where $\nu=\frac{2t-t_0+r}{t_0+r}$ when $\zeta\leq t-s-2t_0$ and $\nu=-\frac{\zeta-\ntau+3(t_0-t)}{t_0+r}$ when $t-s-2t_0\leq \zeta< \ntau-3(t_0-t)$. 
 
\begin{remark} The MAP-estimate $U_\delta$ takes values in the   Cameron-Martin space of the prior.  
The Cameron-Martin space is the intersection of all
linear subspaces where the random variable $U$  belongs with probability one, and since there may be
 uncountable many such linear spaces, the  Cameron-Martin space may be a zero measurable subset of the space where $U$ takes values, see  \cite{Bogachev1998}.
 In the above settings where $U\sim N(0,C_U)$ and $C_U$ is of order $-2r$ the random variable $U$ takes values in $H^{\ntau}$,  
$\ntau=r-s$
where $s>d/2$, and the Cameron-Martin space containing the MAP-estimate is $H^r$. In the Bayesian setting it is natural that the MAP-estimate can not converge in a smaller space than the one $U$ takes values. However, the same behaviour can be seen also in the deterministic setting when the unknown is in $H^r$ and the MAP-estimate is thought as a Tikhonov regularised solution, see \cite{Kekkonen2014}.
\end{remark}

\begin{example}  
Let us study a simple example in two dimensional torus $\T^2$. We consider a problem 
\begin{align*}  
M_\delta=(I-\Delta)^{-1}U+\delta\Noise
\end{align*}
where $\Noise$ is normalised white Gaussian noise that takes values in $H^{-s}(\T^2)$, $s>1$ and $\delta>0$ is the noise amplitude. The model operator $A=(I-\Delta)^{-1}$ is elliptic operator, smoothing of order $2$.  Let us consider the case when $U$ has a priori distribution $N(0,C_U)$ where $C_U=(I-\Delta)^{-2}$, that is, $r=2$. Then $U$ takes values in $H^{\ntau}(\T^2)$, where $\ntau=r-s<1$, almost surely and $U_\delta\in H^2$. Theorem \ref{theorem:main} guarantees us convergence rate $C\delta$ when $\zeta\leq-3$ and (\ref{converge_2}) when $-3\leq \zeta\leq1$. For example we get the following convergence in $L^2(\T^2)$
\begin{align*}
\E\|U_\delta(\omega)-U(\omega)\|_{L^2(\T^2)}\leq C\delta^{\frac{1}{4}-\epsilon}
\end{align*}
with $\epsilon>0$ arbitrarily small. 
\end{example}

\subsubsection{Convergence results in frequentists setting}\label{sec: frequentist}
In the frequentist case one is often interested in the model 
 \begin{align}\label{frequentist}
M_\delta^\dagger(\omega) = A(u^\dagger)+\delta\Noise(\omega)
\end{align}
where the data $M_\delta^\dagger$ is generated by a `true' solution $u^\dagger\in H^{\ntau}(N)$. Above $\Noise$ is normalised white Gaussian noise and $\delta\in \R_+$ is the noise level. In (\ref{frequentist}) all the randomness of the $M_\delta^\dagger$ comes from the randomness of the noise $\Noise$. We denote  
\beq\label{def udelta dagger}
U_\delta^\dagger(\omega)=T_\delta(M_\delta^\dagger(\omega)).
\eeq     
{\mmattitext In the  frequentist setting we consider the case  where $u^\dagger$  is an arbitrary element of the space $H^{\ntau}(N)$ 
where the random variable $U$ takes values, instead of considering almost every element.
Note that even though a set is measure-theoretically large, that is, has probability 1, it can sill be topologically small (meager, or a set of first Baire category). For a discussion on these issues see \cite{Ghosal1997}.}

We can then study the  convergence $\E_{u^\dagger}\|U_\delta^\dagger(\omega)-u^\dagger\|_{H^{\zeta}}$ or in more frequentist spirit $\E_{u^\dagger}\|U_\delta^\dagger(\omega)-u^\dagger\|_{H^{\zeta}}^2$ where we use notation
\begin{align}\label{frequentis_mean matti}
\E_{u^\dagger}
F(U_\delta^\dagger(\omega),u^\dagger,\omega)=\int_{\Omega}
F(T_\delta(Au^\dagger+\delta\Noise(\omega)) ,u^\dagger,\omega)
%
d\mathbb{P}(\omega), 
\end{align}
that is,
the expectation $\E_{u^\dagger}$ is taken  with respect to the noise $\Noise(\omega)$ and the other terms depending on $\omega$, and $u^\dagger$ is considered as a fixed parameter, c.f.\  \eqref{frequentis_mean}. 
This means that after computing the estimator $U_\delta^\dagger$ using Bayesian  methods we will consider the convergence of the estimator $U_\delta^\dagger$ to the `true' solution $u^\dagger$ which is not thought to be a random draw from the prior any more.


\begin{remark}
Next we consider the frequentist case when in addition it is assumed that  $u^\dagger\in H^{\ntau}(N)$ where $\ntau\geq0$. Note that then $r\geq\frac{d}{2}$. The mean integrated squared error (MISE) of an estimator $\widehat{U}$ is defined 
\begin{align}\label{MISE}
R(\widehat{U},u^\dagger)=\E_{u^\dagger}\|\widehat{U}-u^\dagger\|_{L^2}^2.
\end{align}
The minimax risk $r_\delta(H^{\ntau}(N))$ on the Sobolev space $H^{\ntau}(N)$ is then given by 
\begin{align*}
r_\delta(H^{\ntau}(N))=\inf_{\widehat{U}}\sup_{u^\dagger\in H^{\ntau}(N)}R(\widehat{U},u^\dagger)
\end{align*}
where the infimum is taken over all estimators of the form $\widehat{U}=g(M_\delta^\dagger)$ where $g\in \mathcal{B}(H^{-s},H^{\ntau})$. Here $\mathcal{B}(H^{-s},H^{\ntau})$ is the set of Borel measurable functions from $H^{-s}$ to 
$H^{\ntau}$.
\end{remark}

\begin{theorem}\label{thm:frequentist1} 
Let $r> s>d/2$ and $N$ be a $d$-dimensional closed manifold. Let $u^\dagger\in H^{\ntau}(N)$ where $\ntau=r-s>0$. Assume that $C_U$, the covariance operator of the Gaussian prior, is a self-adjoint, injective and elliptic pseudodifferential operator of order $-2r$.
Let $\Noise(y,\omega)$ be white Gaussian noise on $N$. Consider the measurement 
 \begin{align*}
M_\delta^\dagger(y,\omega) = A(u^\dagger(\cdot))+\delta\Noise(y,\omega),\quad\quad \omega\in\Omega,
\end{align*}
where $A\in H\Psi^{-t,-t_0}$, $t>\max\{0, -\tau+s\}$ and $t\leq t_0\leq t+\ntau/3$, is a hypoelliptic pseudodifferential operator on the manifold $N$ and $A:L^2(N)\to L^2(N)$ is injective. Above $\delta\in \R_+$ is the noise level and $\Noise$ takes values in $H^{-s}(N)$ with some $s>d/2$. 
 
Then there is  $C>0$, independent of $\delta$ and $u^\dagger$, such  that 
\begin{align}\label{thm 2 estimate}
\E_{u^\dagger}\|U_\delta^\dagger(\omega)-u^\dagger\|_{L^2}^2 & \leq C(1+\|u^\dagger\|_{H^{\ntau}}^2)\delta^{\frac{2(\tau-3(t_0-t))}{t_0+r}}. 
\end{align}
\end{theorem}
\medskip
 
Note that the assumptions on $C_U$ in Theorem \ref{thm:frequentist1} imply that $C_U$ is a covariance operator
of a random variable taking values in $H^{\ntau}(N)$, see Remark \ref{matti remark} below.

In the elliptic case $t=t_0$ we can write (\ref{thm 2 estimate}) as
\begin{align}\label{rare_frequentist}
\E_{u^\dagger}\|U_\delta^\dagger(\omega)-u^\dagger\|_{L^2}^2 & \leq C(1+\|u^\dagger\|_{H^{\ntau}}^2)\delta^{\frac{2\ntau}{s+\ntau+t}}.
\end{align}
Since $s=\frac{d}{2}+\epsilon$ the above convergence rate (\ref{rare_frequentist}) agrees, up to $\epsilon>0$ arbitrarily small, with the minimax convergence rate, see \cite{Cavalier2008}. The convergence of confidence regions is considered in Section \ref{sec:crediblesets}.

\section{Generalised random variables}\label{sec:Gaussianexample}  
  
This section is largely based on the work of Lasanen \cite{Lasanen2002, Lasanen2012}; see also Piiroinen \cite{Piiroinen2005}. 

For any  $s\in\R$, let $H^{s}(N)$ be the $L^2$-based Sobolev space equipped with Hilbert space inner product
\begin{eqnarray}\label{eq: inner product in Hs}
( \phi,\psi)_{H^{s}(N)}=\int_{N} ((I-\Delta)^{s/2}\phi)(x)\, ((I-\Delta)^{s/2}\psi)(x)\,dx.
\end{eqnarray}
We also define a dual pairing between $H^{-s}(N)$ and $H^{s}(N)$ 
\begin{eqnarray}\label{eq: dual pairing in Hs}
\langle \phi,\psi \rangle_{H^{-s}(N)\times H^{s}(N)} = \int_{N} \phi(x)\, \psi(x) \,dx
\end{eqnarray}
when $\phi,\psi \in C^\infty_0(N)$.
Note that $H^0(N)=L^2(N)$. We often denote $H^s=H^s(N)$ and $L^2=L^2(N)$.

A generalised Gaussian random variable $V$ takes values in the space of generalised functions, and the pairing $\bra  
V,\phi\cet$ with any test function $\phi\in \mathcal{D}=C^\infty(N)$ is a Gaussian random variable taking values in $\R^d$, see \cite{Rozanov1982}. {\mntext The generalised Gaussian random variables we will consider below are assumed to
take values in some Hilbert space, typically in a Sobolev space $H^s(N)$, where the smoothness
index $s\in \R$ may also be negative. 
Now,} if $V$
takes values in $H^s$ we say that $V$ has the covariance operator $B_V:H^s\to H^s$ if
\begin{equation}\label{covop}
  \expec \big(( V-\expec V,\phi)_{H^s}\, ( V-\expec V,\psi )_{H^s}\big)= ( B_V\phi,\psi)_{H^s},
\end{equation}
with any $\phi,\psi\in H^s$, see \cite{Rozanov1982}. 
We can also define covariance operator to be a mapping $C_V:H^{-s}\to H^s$ 
\begin{equation}\label{covop}
  \expec \Big(\langle V-\expec V,\phi\rangle_{H^s\times H^{-s}}\, \langle V-\expec V,\psi\rangle_{H^s\times H^{-s}}\Big)= \langle C_V\phi,\psi\rangle_{H^s\times H^{-s}},
\end{equation}
with any $\phi,\psi\in H^{-s}$, see \cite{Bogachev1998}. 
The connection between $B_V$ and $C_V$ is 
\ba 
B_V = C_V(I-\Delta)^s : H^s\to H^s.
\ea 
Next we will take a closer look to the generalised white Gaussian noise and introduce the 'white noise paradox' by a simple example.  

\subsection{White noise paradox}\label{sec:whitenoise}  
White noise $\Noise$ can be considered as a measurable map $\Noise:\Omega\to \mathcal D^\prime(N)$ where $\Omega$ is the probability space. Then normalised white noise is a random generalised function $\Noise(y,\omega)$ on $N$ for which the pairings
$\bra  \Noise,\phi\cet_{\mathcal D'\times\mathcal D} $ are Gaussian random variables for all test functions $\phi\in \mathcal D=C^\infty(N)$,
 $\E \Noise=0$, and  
\begin{equation}\label{gauss infinite}
\E\bigg(\bra \Noise,\phi\cet_{\mathcal D'\times\mathcal D} \bra \Noise,\psi\cet_{\mathcal D'\times\mathcal D}\bigg)=\bra I\phi, \psi\cet_{\mathcal D'\times\mathcal D}\quad
\hbox{for }\phi,\psi\in \mathcal{D}.
\end{equation}
We will denote this by $\Noise\sim N(0,I)$. A realisation of $\Noise$ is the generalised function $\varepsilon=\Noise(\,\cdot\,,\omega_0)$  on $N$ with a fixed $\omega_0\in \Omega$.

The probability density
function of white noise $\Noise$ is often {\it formally} written in the form  
\begin{equation}\label{formaldensity}
\pi_\Noise(\varepsilon)\underset{\mbox{\tiny \em formally}}=c\exp\bigg(-\frac{1}{2}\|\varepsilon\|_{L^2(N)}^2\bigg).
\end{equation} 
However, the realisations of the  white Gaussian noise are almost surely not in $L^2(N)$. 
This brings us back to the problem in formula (\ref{contTikh0}). 

\bigskip

\begin{example}
Let $\Noise$ be normalised white Gaussian noise defined on the $d$-dimensional  flat torus 
$\T^d=(\R/(2\pi \Z))^d$. Let $e_{\vec{\ell}}\in L^2(\T^d)$,  $\vec{\ell}=(\ell_1,\ell_2,\dots,\ell_d) \linebreak \in \mathbb{Z}^d$ be an
orthonormal basis of $L^2(\T^d)$ consisting of eigenfunctions of Laplacian, numbered so that $-\Delta e_{\vec{\ell}}=|\vec\ell|^2 e_{\vec{\ell}}$.
Such functions  $e_{\vec{\ell}}(x)$ can be chosen to be
normalised products of the sine and cosine functions $\sin(\ell_jx_j)$ and $\cos(\ell_jx_j)$ that form the standard Fourier basis
of $L^2(\T^d)$. 
The Fourier coefficients of $\Noise$ with respect to this basis are  independent, normally distributed $\R$-valued random variables with variance one, that is, $ \bra \Noise,e_{\vec{\ell}}\cet\sim N(0,1)$.
Then
\ba
\E\|\Noise\|_{L^2(\T^d)}^2=\sum_{\vec{\ell}\in \mathbb{Z}^d}\E|\bra \Noise,e_{\vec{\ell}}\cet|^2=\sum_{\vec{\ell}\in \mathbb{Z}^d} 1=\infty.
\ea
This implies that realisations of $\Noise$ are in $L^2(\T^d)$ with probability zero. However, when $s>d/2$ 
\begin{equation}\label{sum_example1}
\E\|\Noise\|_{H^{-s}(\T^d)}^2=\sum_{\vec{k}\in \mathbb{Z}^d}(1+|\vec{\ell}|^2)^{-s}\E|\bra \Noise,e_\ell\cet|^2<\infty
\end{equation}
and hence $\Noise$ takes values in $H^{-s}(\T^d)$ 
{\newtext almost surely (that is, with probability one)}

On the other hand \cite[Theorem 2]{rozanov1971} implies that if $\|\Noise\|_{H^{-s}(\T^d)}^2<\infty$ almost surely then $\E\|\Noise\|_{H^{-s}(\T^d)}^2<\infty$ which yields $s>d/2$. This concludes that the realisations of white noise $\Noise$ are almost surely in the space $ H^{-s}(\T^d)$
if and only if  $s>d/2$. In particular for $s \leq d/2$ the function $x\mapsto \Noise(x,\omega)$ is in $H^{-s}(\T^d)$ only when $\omega\in\Omega_0\subset\Omega$ where $\mathbb{P}(\Omega_0)=0$.
\end{example}

\subsection{The smoothness of the prior}  
Consider the {\mntext continuum} measurement model $M_\delta=AU+\Noise\delta$ where the operator is now viewed as a smoothing map
$A:H^{\ntau}(N)\to H^{\ntau+t}(N)$ for all $\tau\in\R$.    
We construct the prior by choosing $U$ to be a generalised Gaussian random variable taking values in $H^{\ntau}(N)$ and having expectation $\expec U=0$. First we will, however, give a definition for pseudodifferential operators. 

\begin{definition}\label{def: pseudo_symbol}    
Let $m\in\R$. We define the symbol class $S^m(\R^d,\R^d)$ to consist of $a(x,\xi)\in C^\infty(\R^d, \R^d)$ such that for all multi-indices $\alpha$ and $\beta$ and any compact set $K\subset \R^d$ there is such constant $C_{\alpha,\beta,K}>0$ that   
\ba 
|\partial^\alpha_\xi\partial^\beta_x a(x,\xi)| \leq C_{\alpha,\beta,K}(1+|\xi|)^{m-|\alpha|}, \quad\quad \xi\in\R^d,\ x\in K.
\ea
\end{definition}

\begin{definition}\label{def: pseudo_operator}
Let $Y:U\to \R^d$ be local coordinates of the manifold $N$.  
A bounded linear operator  $A:\mathcal D^\prime (N)\to \mathcal D^\prime(N)$ is called a pseudodifferential operator if
for any local coordinates $Y:U\to \R^d$, $U\subset N$, there is a symbol $a\in S^{m}(\R^d\times\R^d)$ such that
for  $u\in C^\infty_0(U)$  we have 
\begin{align*}
Au(y_1)=\int_N k_A(y_1,y_2)u(y_2)dV_g(y_2)
\end{align*}
where $k_A|_{N\times N\setminus diag(N)}\in C^\infty(N\times N\setminus diag(N))$ and $diag(N)=\{(y,y)\in N\times N\, |\, y\in N\}$. Also when $Y:U\to V\subset\R^d$ are local $C^\infty$-smooth coordinates $k_A(y_1,y_2)$ is given on $U\times U$ by 
\begin{align*}
k_A(Y^{-1}(x_1),Y^{-1}(x_2))=\int_{ \R^d} e^{i(x_1-x_2)\cdot\xi} a(x_1,\xi)d\xi
\end{align*}
where $x_1,x_2\in V\subset \R^d$ and $a=a_V\in S^m(V,\R^d)$.
In this case we will write 
\[A\in \Psi^{m}(N),\]
and say that in local coordinates $Y:U\to V\subset\R^d$ the operator $A$ has the symbol $a(x,\xi)\in S^{m}(V\times\R^d)$. 
\end{definition}

We assume that the covariance operator $C_U\in \Psi ^{-2r}$, that is $C_U$ is smoothing of order $2r$, self-adjoint and elliptic. With given $r\in\R$ we have to choose $\tau\in\R$ so that $B_U=C_U(I-\Delta)^{\ntau}\in\Psi^{-2(r-\ntau)}$ is a trace class operator. An operator $B_U$ is in the trace class $\mathfrak{S}^1(H^{\ntau})$ if 
\begin{equation}\label{tracecondition}
Tr_{H^{\ntau}\to H^{\ntau}}(B_U)=\sum_{j=1}^\infty |\lambda_j|<\infty
\end{equation}
where $\lambda_j$ are the eigenvalues of the operator $B_U$. Condition $B_U\in\mathfrak{S}^1(H^{\ntau})$ guarantees that $\E(U,U)_{H^{\ntau}}<\infty$. 

Let $\nu_j$ be the eigenvalues of $B_U^{-1}\in \Psi^{2(r-\ntau)}$. Counting the geometric multiplicity of the eigenvalues, we arrange the eigenvalues of $B_U^{-1}$ in ascending order as
\ba 
\nu_1\leq\nu_2\leq\cdots\leq\nu_k\leq\cdots 
\ea
Since $B_U^{-1}$ is a self-adjoint elliptic operator with smooth coefficients Weyl's law for elliptic operators tells us that the number $N(\nu)=\#\{\nu_j\ |\ \nu_j\leq \nu\}$ of the eigenvalues of $B_U^{-1}$ in a closed manifold less than or equal to $\nu$ has asymptotics
\ba 
N(\nu)\sim c\nu^{\frac{d}{2(r-\ntau)}}\big(1+O(\nu^{-\frac{1}{2(r-\ntau)}})\big)\qquad \text{when}\ \ \nu\to \infty.
\ea

Hence for the eigenvalues $\lambda$ of the operator $B_U\in \Psi^{-2(r-\ntau)}$
\ba 
\lambda_j \sim c' j^{-\frac{2(r-\ntau)}{d}}(1+o(1))\qquad \text{when}\ \ j\to \infty.
\ea
To satisfy condition (\ref{tracecondition}) we require 
\ba 
\sum_{j=1}^\infty |\lambda_j|\leq C\sum_{j=1}^\infty j^{-\frac{2(r-\ntau)}{d}}<\infty
\ea
which gives us the condition $\ntau<-(d/2-r)$. From here on we will assume that $\ntau=r-s<r-d/2$.    

\begin{remark}\label{matti remark} Any elliptic operator  $C_U\in \Psi ^{-2r}$ that defines a non-negative symmetric operator
$C_U:  \mathcal D(N)=C^\infty(N)\to \mathcal D^\prime(N)$ has the property that $B_U=C_U(I-\Delta)^{\ntau}$ is in $\mathfrak{S}^1(H^{\ntau}(N))$.
By \cite{Bogachev1998}, these yield that $C_U$ is a covariance operator of a random variable taking values in 
$H^{\ntau}(N)$.
\end{remark}

The operator $B_U$ corresponds formally to the smoothness prior
\begin{eqnarray*}
\pi_{pr}(u)\underset{\mbox{\tiny \em formally}} = c\exp\bigg(-\frac{1}{2}( B_U^{-1}u,u)_{H^{\ntau}}\bigg) = c\exp\bigg(- \frac{1}{2} \|C_U^{-1/2}u\|_{L^2}^2\bigg)
\end{eqnarray*}
Notice that the realisations of $U$ are almost surely not $r$ times differentiable. In a case $r<d/2$ the realisations of $U$ are almost surely not even in $L^2$ let alone differentiable. 
This is why we need to consider $U$ as taking values in some space $H^{\ntau}(N)$ with possible negative smoothness index $\ntau$.

\section{Proof of the main Theorems}\label{sec:generalproof}  

Before we move to prove theorem \ref{theorem:main} we will give a short introduction to hypoelliptic pseudodifferential operators.  

\begin{definition}\label{def: hypoelliptic}
Let $t,t_0\in\R$. We define symbol class $HS^{-t,-t_0}(\R^d,\R^d)$ to consist of $a(x,\xi)\in C^\infty(\R^d, \R^d)$ for which
\begin{enumerate}
 \item For an arbitrary compact set $K\subset \R^d$ 
 we can find such positive constants $R$, $c_1$ and $c_2$ that
\begin{equation*}
c_1(1+|\xi|)^{-t_0}\leq|a(x,\xi)|\leq c_2 (1+|\xi|)^{-t}, \quad\quad |\xi|\geq R,\ x\in K.
\end{equation*}
\item For any compact set $K\subset \R^d$ 
there exist constants $R$ and $C_{\alpha,\beta,K}$ such that for all multi-indices $\alpha$ and $\beta$
\ba 
|\partial^\alpha_\xi\partial^\beta_x a(x,\xi)|\leq C_{\alpha,\beta,K}|a(x,\xi)|(1+|\xi|)^{-|\alpha|}, \quad\quad |\xi|\geq R,\ x\in K.
\ea
\end{enumerate}
We will denote by $ H\Psi^{-t,-t_0}(N)$ the class of $\Psi$DO with local symbol $a(x,\xi)\in HS^{-t,-t_0}(V\times\R^d)$, see Definition \ref{def: pseudo_operator}. 
\end{definition}

We denote $H^r(N)=H^r$ and $L^2(N)=L^2$ where $N$ is a closed manifold and $\dim N=d$. 


The proof of Theorem \ref{theorem:main} is rather long and technical so we will start by going through the main steps of it in a nutshell. 
The approximated solution we are studying is of the form
\begin{equation*}
  T_{\delta}(m_\delta):=\argmin_{u\in H^r(N)}
  \big\{\|A u\|_{L^2(N)}^2 - 2\langle m_\delta,A u \rangle +\delta^2\|C_U^{-1/2}u\|_{L^2(N)}^2 \big\}.
\end{equation*}
As mentioned before the solution to this is 
\begin{equation}\label{regsolutio}
  T_{\delta}(m_\delta)=( A^*A+ \delta^2C_U^{-1})^{-1}A^*m_\delta.
\end{equation}
We can rewrite the above as 
\begin{equation} \label{uparts}  
\begin{split}
T_{\delta}(m_\delta) & = Z_\delta^{-1}A^*Au+Z_\delta^{-1}A^*(\varepsilon\delta)\\
& = u-\delta^2 Z_\delta^{-1}C_U^{-1}u+Z_\delta^{-1}A^*(\varepsilon\delta)
\end{split}
\end{equation}
where $Z_\delta=A^*A+\delta^2C_U^{-1}$.

To study the convergence of the last term on the right hand side of (\ref{uparts}) we would like to write it in the form 
\begin{align}\label{operator F}
Z_\delta^{-1}A^*(\varepsilon\delta)=\delta^{-1}F^{-1}_\delta(A^*A)^{-1}A^*\varepsilon
\end{align}
where  $F_\delta=(A^*A)^{-1}C^{-1}_U+\delta^{-2}$. In order to show that (\ref{operator F}) is well-defined we first need to prove that $A^*A$ and $F_\delta$ are invertible.  

Lastly we study the converge of $\delta^2 Z_\delta^{-1}C_U^{-1}u$ and $\delta^{-1}F^{-1}_\delta(A^*A)^{-1}A^*\varepsilon$ to zero in appropriate Sobolev spaces and show that the latter term is always dominating.


We will start by showing that $A^*A\in H\Psi^{-2t,-2t_0}$ is invertible and $(A^*A)^{-1}\in H\Psi^{2t_0,2t}$.
Define $A^*:L^2(N)\to L^2(N)$ as the adjoint of an operator $A:L^2\to L^2$. 
We assumed in Theorem \ref{theorem:main} that $A:L^2\to L^2$ is one-to-one. Since $A^*A:H^r\to H^{r+2t}$, $r\in \R$, is hypoelliptic \cite[Propositions 5.2 and 5.3]{shubin1987} $A^*Au=0\in C^\infty$ implies that $u\in C^\infty$ and hence $Au\in L^2$. Now we see that if $A^*Au=0$ then 
\ba 
0=( A^*Au,u )_{L^2}=( Au,Au )_{L^2}=\|Au\|_{L^2}^2
\ea
which implies $Au=0$ and furthermore $u=0$. Thus the operator $A^*A:H^r(N)\to H^{r+2t}(N)$ is one-to-one.

To study the mapping $A^*A\in H\Psi^{-2t,-2t_0}$ from some Sobolev space $H^r$, $r\in\R$, we define $L_r=A^*A:H^r\to H^{r+2t}$.
The adjoint of $L_r$ is denoted by $L_r'=(A^*A)':H^{-(r+2t)}\to H^{-r}$. Let $\phi,\psi\in C^\infty$. Then
\ba 
\langle A^*A \phi,\psi \rangle_{H^{r+2t}\times H^{-(r+2t)}} = \langle \phi, A^*A\psi \rangle_{H^{r}\times H^{-r}},
\ea
that is, $L_r'=L_{-(r+2t)}$ and hence the adjoint is one-to-one. Now we can conclude that $L_r(H^r)\subset H^{r+2t}$ is a dense subset. Next we will prove two lemmas that show that the operator $A^*A:\mathcal{D}'\to\mathcal{D}'$ is also surjective.
  
Since $A^*A\in H\Psi^{-2t,-2t_0}$ is a hypoelliptic pseudodifferential operator it has a parametrix $B_1\in H\Psi^{2t_0,2t}$ \cite[Theorem 5.1]{shubin1987}. Hence for any $r_0>0$ we get norm estimates 
\begin{equation}\label{estimates}
\left\{ \begin{array}{ll}
\|A^*Au\|_{H^{r+2t}} \leq C_1\|u\|_{H^{r}} \\
\|u\|_{H^{r}} \leq C_2\|A^*Au\|_{H^{r+2t_0}}+C_3\|u\|_{H^{r-r_0}}
\end{array} \right.
\end{equation}
for all $u\in C^\infty$.
Next we will show that $C_3$ is zero.

\begin{lemma}\label{norms}
Let $A^*A\in H\Psi^{-2t,-2t_0}$ be an injective hypoelliptic pseudodifferential operator. Then we have the following estimates
\begin{equation*}
C_1\|A^*Au\|_{H^{r+2t}} \leq \|u\|_{H^{r}} \leq C_2\|A^*Au\|_{H^{r+2t_0}}
\end{equation*}
\end{lemma}

\noindent\textit{Proof.} We get the first inequality since $A^*A$ is continuous linear operator. If the second estimate in (\ref{estimates}) is not valid with $C_3=0$ for any $C_2>0$ then we can choose a sequence $u_j$ such that $\|u_j\|_{H^{r-r_0}}=1$ and $\|u_j\|_{H^{r}}\geq j\|A^*Au_j\|_{H^{r+2t_0}}$. When $j>2C_2$ then  
\ba 
\begin{split}
\|u_j\|_{H^{r}} & \leq C_2\|A^*Au_j\|_{H^{r+2t_0}}+C_3\|u_j\|_{H^{r-r_0}}\\
& \leq \frac{1}{2}\|u_j\|_{H^{r}}+C_3\|u_j\|_{H^{r-r_0}}.
\end{split}
\ea
This gives us 
\ba  
\|u_j\|_{H^{r}}\leq 2C_3\|u_j\|_{H^{r-r_0}}=2C_3.
\ea
Since $r_0>0$, the embedding $H^{r}\hookrightarrow H^{r-r_0}$ is compact. Now there exists a subsequence $u_{j_\ell}$ and such a $w\in H^{r-r_0}$  that $\lim_{\ell\to\infty} u_{j_\ell}=w$ in $H^{r-r_0}$. We assumed that $\|u_{j_\ell}\|_{r-r_0}=1$ which implies $\|w\|_{r-r_0}=1$. On the other hand 
\ba 
j_\ell\|A^*Au_{j_\ell}\|_{r+2t_0}\leq \|u_{j_\ell}\|_{r}\leq 2C_3,
\ea
that is, $\lim_{\ell\to\infty}\|A^*Au_{j_\ell}\|_{r+2t_0}=0$ and because $-r_0+2t\leq 2t_0$
\ba 
\lim_{\ell\to\infty}\|A^*Au_{j_\ell}\|_{r-r_0+2t}=0.
\ea
Since $u_{j_\ell}\to w$ in $H^{r-r_0}$ we also have $A^*Au_{j_\ell}\to A^*Aw$ in $H^{r-r_0+2t}$. Combining the above results we see that $\|A^*Aw\|_{H^{r-r_0+2t}}=0$. Operator $A^*A$ is one to one and hence $w=0$. This is a contradiction since $\|w\|_{H^{r-r_0}}=1$. \qed 

\begin{lemma}\label{subsets}
Let $A^*A: H^{r}\to H^{r+2t}$ be an injective hypoelliptic operator. Then the image of $H^{r}$ in the map $A^*A$ satisfies
\begin{equation*}  
H^{r+2t_0} \subset A^*A(H^{r}) \subset H^{r+2t}.
\end{equation*}
\end{lemma}

\noindent\textit{Proof.} The second inclusion is a direct consequence of the mapping properties of $A^*A$.
Let $f\in H^{r+2t_0}$. Since $C^\infty\subset H^{r+2t_0}$ is a dense subset we can find such a sequence $f_j\in C^\infty$ that $\lim_{j\to\infty}f_j=f$ in $H^{r+2t_0}$. Since $A^*A(H^r)\subset H^{r+2t}$ is dense we can also choose a sequence $h_{j,\ell}=A^*Ag_{j,\ell}\in A^*A(H^{r})$ such that  $\lim_{\ell\to\infty}h_{j,\ell}=f_j$ in $H^{r+2t_0}$. Denote $g_j=g_{j,\ell_j}\in H^{r}$ for which $\lim_{j\to\infty}A^*Ag_j=f$ in $H^{r+2t_0}$. Using Lemma \ref{norms} we see 
\ba 
\lim_{j,k\to \infty}\|g_j-g_k\|_{H^r}\leq C_2\lim_{j,k\to \infty}\|A^*Ag_j-A^*Ag_k\|_{H^{r+2t_0}}=0.
\ea
Hence also $g_j\in H^r$ is a Cauchy sequence. Thus there exists such $g\in H^{r}$ that $\lim_{j\to\infty}g_{j}=g$ in $H^{r}$.  
On the other hand, 
\ba 
\lim_{j\to \infty}\|A^*Ag_{j}-A^*Ag\|_{H^{r+2t}}\leq C_1\lim_{k\to \infty}\|g_{j}-g\|_{H^{r}}=0.
\ea
Combining the above we get $A^*Ag=f$. \qed\\ 
 
Using Lemma \ref{subsets} we see that $A^*A(\mathcal{D'})=\mathcal{D'}$, that is the operator $A^*A$ is also onto. 
Now we can conclude that there exists an inverse operator $(A^*A)^{-1}:\mathcal{D'}\to \mathcal{D'}$. It remains to show that the inverse operator is a hypoelliptic pseudodifferential operator.

\begin{lemma} 
A self-adjoint, smoothing, one-to-one hypoelliptic operator $A^*A\in H\Psi^{-2t,-2t_0}$ has an inverse operator $(A^*A)^{-1}\in H\Psi^{2t_0,2t}$.
\end{lemma}

\noindent\textit{Proof.} Denote $B=(A^*A)^{-1}:\mathcal{D'}\to \mathcal{D'}$.
For an operator $A^*A:H^{r}\to H^{r+2t}$ we define $B_0\subset B$ with domain
\begin{equation*}
\mathcal{D}(B_0) =\{f\in H^{r+2t}\ |\ Bf\in H^{r}\} = A^*A(H^{r}).
\end{equation*} 
Using the hypoellipticity of $A^*A$ we see that $A^*Au=f\in C^\infty$ implies $u\in C^\infty$. This gives us $B_0:C^\infty\to C^\infty$. Since $C^\infty$ is a Frech\'et space and $A^*A$ is continuous and linear $A^*A:C^\infty \to C^\infty$ is an open mapping \cite[Theorem 2.11]{Rudin1987}. 
Hence the operator $B_0:C^\infty \to C^\infty$ is continuous.

Since $A^*A$ is hypoelliptic it has a parametrix $B_1\in H\Psi^{2t_0,2t}$ \cite[Theorem 5.1]{shubin1987} 
\begin{equation*}
 \left\{ \begin{array}{ll}
B_1(A^*A) = I+K_1, & K_1\in \Psi^{-\infty}\\
(A^*A)B_1 = I+K_2, & K_2\in \Psi^{-\infty}
\end{array} \right.
\end{equation*}
and we can write
\begin{equation*}
B_0 = B_0((A^*A)B_1-K_2)=B_1-B_0K_2.
\end{equation*}
The operator $B_0K_2:\mathcal{D'}\to C^\infty$ is continuous  and thus we have shown that 
\begin{equation*}
B_0 = B_1  \mod \Psi^{-\infty}.
\end{equation*}  
That is $B_0\in H\Psi^{2t_0,2t}$.
\qed\\  


Next we will examine $\Psi$DOs that depend on spectral variable $\lambda=\delta^{-2}$. For the general theory see \cite{shubin1987}. 

\begin{definition}\label{def: spectral}  
The symbol class $S^m_p(\R^d\times\R^d,\R_+)$ consist of the functions $a(x,\xi,\lambda)$ such that  
\begin{enumerate}
  \item $a(x,\xi,\lambda_0)\in C^\infty(\R^d\times\R^d)$ for every fixed $\lambda_0\geq0$ and
  \item for arbitrary multi-indices $\alpha$ and $\beta$ and any compact set $K\subset \R^d$ there exists a constant 
  $C_{\alpha,\beta,K}$ such that 
  \[|\partial^\alpha_\xi\partial^\beta_x a(x,\xi,\lambda)|\leq C_{\alpha,\beta,K}(1+|\xi|+|\lambda|^{1/p})^{m-|\alpha|}\]  
  for $x\in K$, $\xi\in \R^d$ and $\lambda\geq0$.
\end{enumerate}
We denote by $\Psi_p^m(N,\R_+)$ the class of pseudodifferential operators $A_\lambda$ for which the local symbol $a(x,\xi,\lambda)\in S^m_p(V\times\R^d,\R_+)$, see Definition \ref{def: pseudo_operator}.
\end{definition}

\begin{definition}\label{def: spectral_hypo}
If there are constants $C_1,C_2,R>0$  such that the symbol $a(x,\xi,\lambda)\in S^m_p(\R^d\times \R^d,\R_+)$ satisfies
\ba 
C_1(|\xi|+|\lambda|^{1/p})^{m_0}\leq |a(x,\xi,\lambda)|\leq C_2(|\xi|+|\lambda|^{1/p})^m,
\ea 
for $|\xi|+|\lambda|\geq R$, we say that $a$ is hypoelliptic with parameter $\lambda$ and denote $a(x,\xi,\lambda)\in HS^{m,m_0}_p(\R^d\times \R^d,\R_+)$. We will denote by $H\Psi^{m,m_0}_p(N,\R_+)$ the class of $\Psi$DOs depending on the parameter $\lambda$ whose local symbol belongs to $HS^{m,m_0}_p(V\times \R^d,\R_+)$, see Definition \ref{def: pseudo_operator}.
\end{definition}

Next we will prove that 
\ba   
F_\lambda=( A^*A)^{-1}C_U^{-1}+\lambda 
\ea      
is invertible. Operator $F_\lambda\in H\Psi^{2(t_0+r),2(t+r)}(N)$ is hypoelliptic since $( A^*A)^{-1}C_U^{-1}\in H \Psi^{2(t_0+r),2(t+r)}(N)$ is hypoelliptic and $\lambda\ I\in \Psi^0(N)$. Denote $Q=( A^*A)^{-1}C_U^{-1}$ and its symbol $q(x,\xi)\in HS^{2(t_0+r),2(t+r)}(N)$. Then for the symbol $\sigma(F_\lambda)(x,\xi)=q(x,\xi)+\lambda$ of the operator $F_\lambda$ 
\ba 
|\partial^\a_\xi \partial^\b_x (q(x,\xi)+\lambda)| \leq C_{\a,\b}(1+|\xi|+|\lambda|^{1/(2(t_0+r))})^{2(t_0+r)-|\a|}.
\ea
By \cite[Theorem 9.2.]{shubin1987}  there exist $R>0$ such that for $|\lambda|\in [R,\infty)$ the operator $F_\lambda\in H\Psi^{2(t_0+r),2(t+r)}_{2(t_0+r)}(N,\R_+)$ is invertible with
\ba 
F^{-1}_\lambda\in H\Psi^{-2(t+r),-2(t_0+r)}_{2(t_0+r)}(N,[R,\infty)).
\ea

We have now shown that the operator $Z_\delta$ can be written    
\begin{equation}\label{def Z}
Z_\delta^{-1} = \lambda\bigg(( A^*A)^{-1}C_U^{-1}+\lambda I\bigg)^{-1}( A^*A)^{-1}
\end{equation}
where $\lambda=\delta^{-2}$.
Hence we can rewrite (\ref{uparts})   
\begin{align}\label{uparts2}
T_\lambda(m) & = u-\lambda^{-1} Z_\delta^{-1}C_U^{-1}u+\sqrt{\lambda} F_\lambda^{-1}(A^*A)^{-1}A^*\varepsilon.
\end{align}

Now we will proceed to study the convergence of the second and third term on the right hand side of (\ref{uparts2}).
For the third term of (\ref{uparts2}) we have  $(A^*A)^{-1}A^*:H^{-s}\to H^{k}$, $k=-s+t-2t_0$ and $F_\lambda^{-1}:H^{k}\to H^{k+2(t+r)}$. Hence when $\zeta\leq k+2(t+r)$ we  have
\ba 
\|F_\lambda^{-1}(A^*A)^{-1}A^*\varepsilon\|_{H^{\zeta}} 
&\leq & \| F_\lambda^{-1}\|_{k,\zeta} \|(A^*A)^{-1}A^*\varepsilon\|_{H^{k}}
\ea  
where $\| F^{-1}_\lambda\|_{k,\zeta}$ is the norm of $F^{-1}_\lambda:H^{k}(N)\to H^{\zeta}(N)$ and $k,\zeta\in \R$. 
Next we want to study what happens to the norm when $\lambda\to\infty$. 

We have the following norm estimates for  $F^{-1}_\lambda\in \Psi^{m}_p(N,\R_+)$ when $ \ell\geq m$ and $\lambda$ large enough \cite[Theorem 9.1.]{shubin1987}
\begin{align}
\| F^{-1}_\lambda\|_{k,k-\ell}\leq C_{k,\ell}(1+|\lambda|^{1/p})^m, \quad\quad \text{if}\quad \ell\geq 0 \label{Shubin1}\\
\| F^{-1}_\lambda\|_{k,k-\ell}\leq C_{k,\ell}(1+|\lambda|^{1/p})^{-(\ell-m)}, \quad\quad \text{if}\quad \ell\leq 0. \label{Shubin2}
\end{align} 
In our case $ F_\lambda^{-1}\in\Psi^m_p$ where $m=-2(t+r)$ and $p=2(t_0+r)$. We will write $\| F_\lambda^{-1}\|_{k,\zeta}=\| F_\lambda^{-1}\|_{k,k-\ell}$ where $\ell=k-\zeta\geq m$.

First we study the case when $\ell\geq0$ that is $\zeta\leq k$. Inequality (\ref{Shubin1}) gives us the norm estimate 
\ba 
\|F^{-1}_\lambda\|_{k,k-\ell}\leq C(1+|\lambda|^{1/p})^m.
\ea 
Because we want $\sqrt{\lambda}\|F_\lambda^{-1}(A^*A)^{-1}A^*\varepsilon\|_{H^{\zeta}}$ to converge when $\lambda\to \infty$ we have to require that 
\ba
\frac{m}{p}=\frac{-2(t+r)}{2(t_0+r)}< -\frac{1}{2}.
\ea
This is true when $t_0<2t+r$. 
 
When $\ell\leq0$ we have $k\leq \zeta \leq k+2(t+r)$ and can use (\ref{Shubin2})
\ba
\|F^{-1}_\lambda\|_{k,k-\ell}\leq C(1+|\lambda|^{1/p})^{-(\ell-m)}.
\ea
 For convergence we need 
\ba 
\frac{m-\ell}{p}=\frac{-2(t+r)-k+\zeta}{2(t_0+r)}< -\frac{1}{2}
\ea
that is $k\leq \zeta<k+2t+r-t_0=r-s-3(t_0-t)$ which can be true only if $t_0< 2t+r$.

Next we will prove the convergence of the term ${\delta^2} Z_\delta^{-1}C_U^{-1}u$ in $H^{\zeta}$.
Since we got above that $\zeta<\ntau-3(t_0-t)$ we can write $\zeta=\ntau-\theta$ where $\theta\geq 3(t_0-t)\geq 0$.
We need to find such $\eta\geq0$ and $\gamma\geq0$ that $\gamma+\eta=1$ and $t_0\gamma-r\eta+r-\theta/2=0$. Define $\gamma=\theta/2(t_0+r)$ and $\eta=1-\theta/2(t_0+r)$. Now $\eta\geq 0$ only if $\theta\leq 2(t_0+r)$. Hence we will choose $\theta=\min\{\ntau-\zeta,2(t_0+r)\}$. 

Since $Z_\delta=A^*A+ {\delta^2}C_U^{-1}$ where $A^*A\geq c_1(I-\Delta)^{-t_0}$ and $c_2(I-\Delta)^{r}\leq C_U^{-1} \leq c_3(I-\Delta)^{r}$ we get
\beq\label{Matti 2}
\begin{split}
\|{\delta^2} Z_\delta^{-1}C_U^{-1}u\|_{H^{\zeta}} 
& \leq \delta^2\|(A^*A)^{-\gamma}(\delta^2 C_U^{-1})^{-\eta} c_3(I-\Delta)^{r+\frac{\zeta}{2}} u\|_{L^2}\\
& \leq \delta^2\|(c_1(I-\Delta)^{-t_0})^{-\gamma}(c_2\delta^2 (I-\Delta)^{r})^{-\eta} (I-\Delta)^{r-\frac{\theta}{2}+\frac{\ntau}{2}} u\|_{L^2}\\
& = C \delta^{\frac{\theta}{t_0+r}} \| u\|_{H^{\ntau}}
\end{split}
\eeq
where $\theta=\min\{\ntau-\zeta,2(t_0+r)\}$.

%

Adding the above results together we can prove Theorem \ref{theorem:main}}.  
\medskip

\noindent\textit{Proof of Theorem \ref{theorem:main}.} To get the speed of convergence we use the fact that $U$ and $\Noise$ are independent. 
Similarly to (\ref{uparts2}) we get 
\beq\label{uparts2 Matti}
U_\delta(\omega)-U(\omega)=
 -\delta^2 Z_\delta^{-1}C_U^{-1} U(\omega)+\frac 1\delta F_{\delta^{-2}}^{-1}(A^*A)^{-1}A^*\Noise(\omega),\hspace{-1cm}
\eeq
where by (\ref{Matti 2}),
\beq\label{noise Matti1}
\|\delta^2 Z_\delta^{-1}C_U^{-1} U\|_{H^{\zeta}}\hspace{-3mm} &\leq&\hspace{-3mm} C \delta^{\frac{\theta}{t_0+r}} \| U\|_{H^{\ntau}}.
\eeq
For the second part on the right hand side of \eqref{uparts2 Matti} we can write
\beq
\label{noise Matti2}
\E \|\frac 1\delta F_{\delta^{-2}}^{-1}(A^*A)^{-1}A^*\Noise(\omega)\|^p_{H^{\zeta}}\hspace{-3mm}  &\leq&\hspace{-3mm} 
\|\frac 1\delta  F_{\delta^{-2}}^{-1}(A^*A)^{-1}A^*\|_{H^{-s}\to H^{\zeta}}^p  \,\E \|   \Noise(\omega)\|^p_{H^{-s}}, \hspace{-1.9cm}
\eeq
with $p\in \{1,2\}$ and $\theta=\min\{\ntau-\zeta,2(t_0+r)\}$.

When $\zeta\leq t-s-2t_0$ we get
\ba   
\E\|U_\delta(\omega)-U(\omega)\|_{H^{\zeta}} 
&\leq & \delta^2\E\| Z_\delta^{-1}C_U^{-1}U(\omega)\|_{H^{\zeta}}+\delta^{-1}\E\| F_\delta^{-1}(A^*A)^{-1}A^*\Noise(\omega)\|_{H^{\zeta}}\\
&\leq & C_1\delta^{\frac{\theta}{t_0+r}}\E\|U(\omega)\|_{H^{\ntau}}+C_2\delta^{-1+\frac{2(t+r)}{t_0+r}}\E\|\Noise\|_{H^{-s}}\\
&\leq & C\max\Big\{\delta^{\frac{\theta}{t_0+r}},\delta^{\frac{2t+r-t_0}{t_0+r}}\Big\}\\  
\ea
where $\theta=\min\{\ntau-\zeta,2(t_0+r)\}$. Next we will study which of the terms is dominating. The noise term $\delta^{-1}F_\delta^{-1}(A^*A)^{-1}A^*\Noise(\omega)$ is dominating if 
\begin{align*}
2t+r-t_0\leq \theta.
\end{align*}
Assume first that $\theta=\ntau-\zeta$. Then 
\begin{align*}
\theta=\ntau-\zeta\geq 2t_0+r-t\geq 2t+r-t_0.
\end{align*}
If $\theta=2(t_0+r)$ we get
\begin{align*}
\theta=2(t_0+r)\geq 2t_0+t_0-2t+r\geq 2t+r-t_0
\end{align*}
since $t\leq t_0< 2t+r$.
Hence the noise term is dominating in both cases and we have proven
\ba 
\E\|U_\delta(\omega)-U(\omega)\|_{H^{\zeta}} \leq C\delta^{\frac{2t-t_0+r}{t_0+r}}.
\ea

If $t-s-2t_0\leq \zeta< \ntau-3(t_0-t)$ we get 
\ba  
\E\|U_\delta(\omega)-U(\omega)\|_{H^{\zeta}} 
&\leq & \delta^2\E\| Z_\delta^{-1}C_U^{-1}U(\omega)\|_{H^{\zeta}}+\delta^{-1}\E\| F_\delta^{-1}(A^*A)^{-1}A^*\Noise(\omega)\|_{H^{\zeta}}\\
&\leq & C_1\delta^{\frac{\ntau-\zeta}{t_0+r}}\E\|U(\omega)\|_{H^{\ntau}}+C_2\delta^{-1+\frac{2(t+r)+\ell}{t_0+r}}\E\|\Noise(\omega)\|_{H^{-s}}\\
&\leq & C\max\Big\{\delta^{\frac{\ntau-\zeta}{t_0+r}},\delta^{\frac{2t+r-t_0+\ell}{t_0+r}}\Big\}\\  
\ea
Above $\ell=t-s-2t_0-\zeta$. Note that when $\zeta\geq t-s-2t_0$ then $\ntau-\zeta\leq 2(r+t_0)$. The noise term is dominating if 
\begin{align*}
2t+r-t_0+\ell\leq \ntau-\zeta.
\end{align*}
This is always true since $\ntau=r-s$ and
\begin{align*}
2t+r-t_0+\ell=r-s-\zeta-3(t_0-t)\leq \ntau-\zeta.
\end{align*}
Hence we can conclude  
\ba  
\E\|U_\delta(\omega)-U(\omega)\|_{H^{\zeta}} \leq  C\delta^{\frac{2t+r-t_0+\ell}{t_0+r}}.
\ea
\qed 
\medskip

\noindent\textit{Proof of Theorem \ref{thm:frequentist1}.}  
Similarly to (\ref{uparts2 Matti}) we get
\beq\label{uparts2 Matti FRE}  
U^\dagger_ \delta(\omega)-u^\dagger=
 -\delta^2 Z_\delta^{-1}C_U^{-1} u^\dagger+\frac 1\delta F_{\delta^{-2}}^{-1}(A^*A)^{-1}A^*\Noise(\omega),\hspace{-1cm}
\eeq
where the first term on the right side satisfies, by (\ref{Matti 2}),
\beq\label{noise Matti1 FRE}
\|\delta^2 Z_\delta^{-1}C_U^{-1} u^\dagger\|_{H^{\zeta}}\hspace{-3mm} &\leq&\hspace{-3mm} C \delta^{\frac{\theta}{t_0+r}} \| u^\dagger\|_{H^{\ntau}}
\eeq
with $\theta=\min\{\ntau-\zeta,2(t_0+r)\}$. The expectation of the second term in the right side is estimated in (\ref{noise Matti2}).
Analysing the obtained terms as in the proof of Theorem \ref{theorem:main}, we obtain
\ba  
\E\|U^\dagger_ \delta(\omega)-u^\dagger\|^2_{H^{\zeta}} \leq   C(1+\|u^\dagger\|_{H^{\ntau}}^2)\delta^{\frac{2(\ntau-3(t_0-t))}{t_0+r}}.
\ea
\qed\\

\section{Posterior distribution and confidence regions}\label{sec:crediblesets}

One advantage Bayesian inversion offers over deterministic regularization is uncertainty quantification. Since the solution to the Bayesian inverse problem is the posterior distribution of the unknown we can study its credible sets and their contraction in some Sobolev space $H^{\zeta}$ when $\delta\to 0$. A Bayesian credible set is a region in the posterior distribution that contains a large fraction of the posterior mass, for instance, 95\%. We are dealing with Gaussian distributions so we define the credible sets to be central regions. This means these sets are defined as central balls with $u_{\delta}$ as a centre. 

The above mentioned credible sets are often used to visualise the remaining Bayesian uncertainty in the estimate. Frequentists use another kind of uncertainty quantification called confidence region. A confidence region is a range of values that frequently includes the unknown of interest if the experiment is repeated. We can define confidence regions as central balls with $u_{\delta}^\dagger$ as the centre. Here $u_{\delta}^\dagger$ is the frequentist approximated solution generated by a true solution $u^\dagger$. How frequently the ball around the approximated solution, with different realisation of the noise, contains the true solution is determined by the confidence level. See for example \cite{Gine2015, Vaart2000}.

In the finite-dimensional parametric case and under mild conditions on the prior Bernstein--von Mises theorem provides that the credible sets of smooth models are asymptotically equivalent with the frequentist confidence regions based on the maximum likelihood estimator, see \cite{Vaart2000}. In infinite-dimensional case there is no corresponding theorem and Bayesian credible sets are not automatically frequentist confidence sets. This means that if we assume that the data is generated by a `true parameter', it is not automatically true that credible sets contain that truth with probability at least the credible level. However the correspondence of Bayesian and frequentist uncertainty has been studied in many recent papers see e.g.  \cite{Castillo2013,Castillo2014,Knapik2011,Leahu2011,Ray2014,szabo2015}. These results are important since they show that some credible sets can give a good idea of the uncertainty of the estimate in the classical sense. In this section we show that the posterior distribution converges and we give some convergence rates. We also prove that in the elliptic case the frequentist posterior contractions rates agrees, up to $\varepsilon>0$ arbitrarily small, with the minimax convergence rate. We do not address the question about the frequentist coverage of the credible sets. 

We will start by studying the convergence of the posterior covariance $C_\delta$ which, with the convergence of the posterior mean $U_\delta$, guarantees the convergence of the posterior distribution. 

When $U\sim N(0,C_U)$, $\Noise\sim N(0,I)$ and 
\begin{align}\label{Bayesian setting}
M_\delta=AU+\delta\Noise 
\end{align}
the conditional probability distribution of $U$ with respect to the measurement $M_\delta$ is a Gaussian measure with mean $U_\delta$ and covariance \cite{Lehtinen1989,Mandelbaum1984} 
\begin{align}\label{Lehtinen}   
C_\delta= C_U-C_UA^*(AC_UA^*+\delta^2I)^{-1}AC_U. 
\end{align}
If $A^*:\mathcal{D}'(N)\to\mathcal{D}'(N)$ is invertible we can rewrite the above
\begin{align}\label{covariance}
\begin{split}
C_\delta & = \delta^2 \big(A^*A+\delta^2C_U^{-1}\big)^{-1}\\
& = \bigg(( A^*A)^{-1}C_U^{-1}+\delta^{-2} I\bigg)^{-1}( A^*A)^{-1}.
\end{split}
\end{align}
Note that the covariance operator is deterministic and thus independent of $M_\delta$.

We define $F_\lambda=( A^*A)^{-1}C_U^{-1}+\lambda I $, where $\lambda=\delta^{-2}$, as in section \ref{sec:generalproof}. Then 
\begin{align*}
F^{-1}_\lambda\in \Psi^{m}_{p}(N,[R,\infty))
\end{align*}
where $m=-2(t+r)$ and $p=2(t_0+r)$.
Using the norm estimate (\ref{Shubin2}) we get
\begin{align*}
\|\delta^2Z_\delta^{-1}\|_{-\ntau,\ntau} & = \|\big(( A^*A)^{-1}C_U^{-1}+\delta^{-2} I\big)^{-1}\|_{-\ntau-2t,\ntau}\|( A^*A)^{-1}\|_{-\ntau,-\ntau-2t}\\
& \leq c\|F_\lambda^{-1}\|_{k,k-\ell}\\
& \leq c(1+\lambda^{\frac{1}{p}})^{-(\ell-m)}.
\end{align*}
Above $k=-\ntau-2t$ and $\ell=-2(\ntau+t)<0$. Since $\ntau=r-s$ we can write
\begin{align*}
\frac{\ell-m}{p} & =\frac{-\ntau+r}{t_0+r}\\
& =\frac{s}{t_0+r}
\end{align*}
and hence we get the following convergence rate for the posterior covariance 
\beq\label{cov op convergence}
\|C_\delta\|_{-\ntau,\ntau} & \leq c\delta^\frac{2s}{t_0+r}\\
& \leq c\delta^\frac{d}{t_0+r}.\nonumber
\eeq

We see that the more smoothing the forward operator $A$ is the worse convergence we get. Note that $r$ and $s$ do not only affect the convergence speed but also the spaces between which the norm is taken.

\begin{remark}  
Observe that the random variable $U$ takes values in $H^{\ntau}$ and the estimate (\ref{cov op convergence}) concerns the mapping properties of the posterior covariance operator $C_\delta$ from the dual space $ H^{-\ntau}=( H^{\ntau})^\prime$ to the space
$ H^{\ntau}$. For strictly positive $\delta>0$ the MAP estimator $U_\delta$  belongs to the space $H^r$, $r=\ntau+s\geq \ntau$,
but as $\delta\to 0$, the  MAP estimators $U_\delta$ converge in a less regular space $H^\zeta$, $\zeta<\ntau-3(t_0-t)\leq \ntau$,
see (\ref{zeta convergence}).
\end{remark}

%
%

\subsection{Contraction of the posterior distribution}
Next we consider the inverse problem using the frequentist setting described in Subsection \ref{sec: frequentist} with the additional assumption that $\ntau>0$. We assume below that  $C_U$ satisfies the assumptions in Theorem \ref{thm:frequentist1} that in particular imply that $C_U$ is the covariance operator
of a random variable $U$ taking values in $H^{\ntau}(N)$, see Remark \ref{matti remark}.
We recall
that we consider a fixed `true'  solution $u^\dagger\in H^{\ntau}(N)$ and the noise model 
$M^\dagger_\delta(\omega)=Au^\dagger+\delta \Noise(\omega)$ as in \eqref{frequentist}. Also, note that the
MAP-estimate is then $U^\dagger_\delta=T_\delta(M^\dagger_\delta)$.

In the frequentist case one is often interested in the the limiting behaviour of the posterior measure ${P}_{M_\delta^\dagger}$  when $\delta\to0$. {\mmattitext Here, ${P}_{M_\delta^\dagger}$ is a random measure  in $H^{\ntau}(N)$, depending on $\delta$ and
the MAP estimator  $U_\delta^\dagger=
U_\delta^\dagger(\omega)=T_\delta(M_\delta^\dagger(\omega))$ (that further depends on the deterministic variable
$u^\dagger$  and
the realisation $\mathcal E(\omega)$ of the random noise).
Let $W_\delta$  be
a Gaussian random variable, taking values in $H^{\ntau}(N)$,
 that  is independent of the noise $\Noise$,
has zero mean    and
the covariance operator $C_\delta$, see (\ref{covop}).
For a measurable set $B\subset H^{\ntau}(N)$ we define
\beq\label{Post measure}
{P}_{M_\delta^\dagger}(B)=\mu_\delta\Big(\Big\{w_\delta\in H^{\ntau}(N)\, |\, w_\delta=b-U_\delta^\dagger(\omega)\hbox{ with }b\in B\Big\}\Big)\hspace{-1cm}
\eeq
where $\mu_\delta=N(0,C_\delta)$. Roughly speaking, ${P}_{M_\delta^\dagger}$ is a Gaussian
measure
in $H^{\ntau}(N)$ with the mean  $U_\delta^\dagger$  and
the covariance operator $C_\delta$. 

Recall that we consider the probability space $(\Omega,\Sigma,\mathbb{P})$ and denote by $\chi_S$
 the indicator function of $S$.
%
Let $S\in \Sigma$. We use the notations
\beq\label{con prob}
\E_{M_\delta^\dagger} (F(\omega,u^\dagger))
:=\E(F(\omega,u^\dagger)| \mathcal F),\quad
\mathbb{P}_{M_\delta^\dagger}(S)
:=\mathbb{P}(S| \mathcal F)=\mathbb{E}(\chi_S| \mathcal F),\hspace{-1cm}
\eeq
for the conditional expectation and conditional probability. Above $\mathcal F\subset \Sigma$ is the $\sigma$-algebra  generated by random variable $M_\delta^\dagger(\omega)$ or equivalently, the noise $\Noise(\omega)$. Roughly speaking, in the notation $\E_{u^\dagger} F$ the subindex $u^\dagger$ reminds that $u^\dagger$ is a fixed parameter and the expectation is taken only with respect the noise. The notation ${P}_{M^\dagger_\delta}(B)$ indicates that the measure of $B$ is computed using the posterior probability measure which mean $U^\dagger_\delta=T_\delta(M^\dagger_\delta)$ depends on the measurement $M^\dagger_\delta$. 
Since the random variable $W_\delta$ has distribution $\mu_\delta$, we have by (\ref{Post measure})
\beq\label{Post measure2}
{{P}}_{M_\delta^\dagger}(B)=\int_{H^\ntau(N)} \chi_B( w_\delta+U_\delta^\dagger)\, d\mu_\delta(w_\delta)
=\mathbb{P}_{M_\delta^\dagger}(\{W_\delta+U_\delta^\dagger\in B\}).
\hspace{-1cm}\eeq

Following the approach in \cite{Florens2016,Ghosal2000,Knapik2011,Vaart2008} we next show that the posterior measure contracts to a Dirac measure centred on the fixed true
solution $u^\dagger$.




\begin{theorem}\label{Thm:contraction}
Let $r> s>d/2$ and $N$ be a $d$-dimensional closed manifold. Let $u^\dagger\in H^{\ntau}(N)$ where $\ntau=r-s>0$. Assume that $C_U$, the covariance operator of the Gaussian prior, is a self-adjoint, injective and elliptic pseudodifferential operator of order $-2r$.
Let $\Noise(y,\omega)$ be white Gaussian noise on $N$. Consider the measurement 
 \begin{align*}
M_\delta^\dagger(y,\omega) = A(u^\dagger(\cdot))+\delta\Noise(y,\omega),\quad\quad \omega\in\Omega,
\end{align*}
where $A\in H\Psi^{-t,-t_0}$, $t>\max\{0, -\tau+s\}$ and $t\leq t_0\leq t+\ntau/3$, is a hypoelliptic pseudodifferential operator on the manifold $N$ and $A:L^2(N)\to L^2(N)$ is injective. We assume also that $A^*:\mathcal{D}'(N)\to\mathcal{D}'(N)$ is invertible. Above $\delta\in \R_+$ is the noise level and $\Noise$ takes values in $H^{-s}(N)$ with some $s>d/2$. 
Let 
$U_\delta^\dagger$ be the MAP estimated given by (\ref{def udelta dagger}).
 
Let $\kappa<\kappa_0=\frac{2(\ntau-3(t_0-t))}{t_0+r}$, $c_0>0$, and $R>0$.
Then  there is $c_1>0$  such that   \begin{align}\label{contraction}
\sup_{\|u^\dagger\|_{H^{\ntau}}\leq R}\E_{u^\dagger}{P}_{M_\delta^\dagger}\Big\{ u\in H^{\ntau}(N)\, \big|\, \|u-u^\dagger\|_{L^2(N)}\geq c_0\delta^{\kappa}\Big\}\leq c_1 \delta^{2(\kappa_0-\kappa)}\to 0,
\end{align}
as $\delta\to 0$.
\end{theorem}

\noindent\textit{Proof.} 
{\mmattitext  Let $u^\dagger\in H^{\ntau}(N)$ and $W_\delta$ be the Gaussian variable defined above.
%
%
%
Using the Markov inequality and (\ref{Post measure2}), we get
\beq\nonumber   
 P_{M_\delta^\dagger}\Big\{ u\in H^{\ntau}(N)\,\big|\, \|u-u^\dagger\|_{L^2}\geq c_0\delta^{\kappa}\Big\}
\hspace{-2mm} &=&\hspace{-2mm} \mathbb{P}_{M_\delta^\dagger}\Big(\Big\{\|W_\delta+U_\delta^\dagger-u^\dagger\|_{L^2(N)}\geq c_0\delta^{\kappa}\Big\}\Big) \hspace{-1cm}\\ \label{PM Markov}
&\leq&\hspace{-2mm}    \frac{1}{(c_0\delta^\kappa)^2} 
\E_{M_\delta^\dagger}\big( \|W_\delta+U_\delta^\dagger- u^\dagger\|_{L^2(N)}^2\big).\hspace{-1cm}
\eeq
Since $W_\delta$ and $U_\delta^\dagger$ are independent and $W_\delta$  has the covariance operator $C_\delta$,
we obtain using notations (\ref{con prob}) 
\beq\nonumber
\E_{u^\dagger} \E_{M_\delta^\dagger}( \|W_\delta+U_\delta^\dagger- u^\dagger\|_{L^2(N)}^2)
\hspace{-2mm} &=&\hspace{-2mm} 
\E_{u^\dagger} \E_{M_\delta^\dagger} \|W_\delta\|_{L^2}^2+\E_{u^\dagger} \E_{M_\delta^\dagger} \|U_\delta^\dagger- u^\dagger\|_{L^2}^2
\\
\hspace{-2mm} &=&\hspace{-2mm} \label{SPC}
Tr_{L^2(N)\to L^2(N)}(C_\delta)+\E_{u^\dagger}\|U_\delta^\dagger-u^\dagger\|_{L^2(N)}^2.
\eeq

  We have shown in Theorem \ref{thm:frequentist1} that the second term on the right side of  \eqref{SPC} can be estimated by $c_2(1+R^2)\delta^{2\kappa_0}$
  with some $c_2>0$. Hence it is enough to show that 
\begin{align*}
Tr_{L^2\to L^2}(C_\delta)\leq c\delta^{2\kappa_0}
\end{align*}
with $\kappa_0=\frac{2(\ntau-3(t_0-t))}{t_0+r}$.} We can estimate the trace by writing 
\begin{align*}
Tr_{L^2\to L^2}(C_\delta) &= Tr_{L^2\to L^2}\big((I-\Delta)^{-s}(I-\Delta)^{s}C_\delta\big)\\
&\leq Tr_{L^2\to L^2}\big((I-\Delta)^{-s}\big)\|(I-\Delta)^{s}C_\delta\|_{L^2\to L^2}\\
\end{align*}
Above $(I-\Delta)^{-s}$ is trace class operator in $L^2$ since 
\ba 
\sum_{j=1}^\infty(j^{2/d})^{-s}<\infty \quad \text{when}\ s>d/2.
\ea

As before we get
\begin{align*}
\|(I-\Delta)^{s}C_\delta\|_{0,0} 
& = \|(I-\Delta)^{s}\|_{2s,0}\|\F_\lambda^{-1}\|_{-2t,2s}\|( A^*A)^{-1}\|_{0,-2t}\\
& \leq c\|F_\lambda^{-1}\|_{-2t,-2t-\ell}\\
& \leq c(1+\lambda^{\frac{1}{p}})^{-(\ell-m)}.
\end{align*}
Above $\ell=-2(s+t)<0$. We can write
\begin{align*}
\frac{\ell-m}{p} & =\frac{\ntau}{t_0+r}\\
\end{align*}
and hence  
\begin{align*}
Tr_{L^2\to L^2}(C_\delta)\leq c\delta^{\frac{2\ntau}{t_0+r}}\leq c\delta^{2\kappa_0}.
\end{align*}
\qed\\

Note that in the elliptic case $t_0=t$ we get contraction 
\begin{align*}
\E_{u^\dagger}{P}_{M_\delta^\dagger}\Big\{ u\in H^{\ntau}(N)\, \big|\, \|u-u^\dagger\|_{L^2(N)}\geq c_0\delta^{\kappa}\Big\}\to 0
\end{align*}
when $\delta\to 0$ for all $c_0>0$  and $\kappa<\frac{2\ntau}{s+\ntau+t}$. Since $s=\frac{d}{2}+\epsilon$ the above convergence rate agrees, up to $\epsilon>0$ arbitrarily small, with the minimax convergence rate.

\begin{remark}
{\mmattitext 
Above we have assumed that  $u^\dagger$ is in $H^{\ntau}$. This correspond to the fact that the random variable
$U$, having the covariance operator $C_U\in \Psi^{-2r}$, takes values in $H^{\ntau}$. 
The $L^2(N)$ norm in the contraction formula (\ref{contraction})  can be considered as a loss function
on $H^{\ntau}(N)$.
Note that the loss function $d(v_1,v_2)=\|v_1-v_2\|_{L^2}$  defines a distance function in the $H^{\ntau}(N)$, but the obtained metric space is
not complete. When the direct map $A$  is the identity map, similar estimates with different loss functions have been studied in a general setting
in \cite{Hoffmann2015}.  However, from the point of view of inverse problems \cite{Hoffmann2015} corresponds to the case when the direct operator and the covariance operator of the prior commute. This differs from the problem analysed in our paper, where covariance operator $C_U$ and the operator $A$ may not commute, and are of quite different type in the sense that $C_U$ is an elliptic operator but  $A$ is  hypoelliptic operator.
The phenomenon that the solution $u^\dagger$ is assumed to be in a smoother space, in our case in 
$H^{\ntau}$, and the convergence of the posterior distribution is analysed using a loss function given by a less strict norm, in our case  $L^2$-norm, 
appears in many frequentist studies, see e.g. Theorems 2.2 and 2.3 and Remark 3.6 in \cite{Agapiou2013}. 
Conditions similar to  the smoothness requirement $u^\dagger\in H^{\ntau}$ 
are also encountered in classical regularisation theory \cite{Engl1996a} 
where this type of conditions are
called source conditions. 
}
\end{remark}
\subsection{Convergence of the posterior distribution in Bayesian settings} 
Next we will proceed to study the contraction of the posterior distribution using Bayesian techniques, the measurement model
$M_\delta=AU+\delta\Noise$ and the MAP estimator $U_\delta=T_\delta(M_\delta)$. Let us write
$$
V_\delta=U_\delta+W_\delta
$$
where $W_\delta\sim N(0,C_\delta)$ is a Gaussian variable having the covariance 
operator $C_\delta$ given in (\ref{covariance}) and 
the zero mean. Random variables $W_\delta$ and $U_\delta$ are assumed to be independent.
Let $\mathcal M_\delta=\sigma(M_\delta)$ be the $\sigma$-algebra
generated by the random variable $M_\delta$.
Then the distribution of the random variable $V_\delta$ is the same as 
the posterior distribution of $U$ with respect to the $\sigma$-algebra $\mathcal M_\delta$.

Let $\nu_\delta$ be the posterior distribution of $U$ with respect to the $\sigma$-algebra $\mathcal{M}_\delta$. Equivalently $\nu_\delta$ is the distribution of $V_\delta$ in the Sobolev space $H^{\zeta_1}(N)$ where $\zeta_1\leq \ntau$. Let $\mu_\delta$ be the distribution of the random variable $W_\delta$ which is independent of $U_\delta=T_\delta(M_\delta)$. Then the conditional expectation of the indicator function  $\chi_{B_{\zeta_1}(U_\delta(\omega),R(\delta))}(U(\omega))$ with respect to $\mathcal{M}_\delta$ is 
\begin{align}\label{credibleset_1}
\begin{split}
\E\big(\chi_{B_{\zeta_1}(U_\delta,R(\delta))} \,|\, \mathcal{M}_\delta\big)(\omega)
& = \mathbb{P}\big(\{U\in B_{\zeta_1}(U_\delta,R(\delta)),R(\delta))\} \,|\, \mathcal{M}_\delta\big)(\omega)\\
& = \mathbb{P}\big(\{W_\delta\in B_{\zeta_1}(0,R(\delta))\}\big)(\omega)\\
& =\mu_\delta\big( B_{\zeta_1}(0,R(\delta))\big)(\omega).
\end{split}
\end{align}
Above $B_{\zeta_1}(0,R(\delta))$ denotes a ball in $H^{\zeta_1}$ of radius $R(\delta)$.

Let $\mathcal{U}$ be the $\sigma$-algebra generated by $U$. By \cite[Theorem 10.2.2]{Dudley1989} there are regular conditional probabilities ${\bf{P}}(K \,|\, \mathcal{M}_\delta)(\omega)$ for all $K\in \mathcal{B}(H^{\zeta_1})$ and $\omega\in\Omega$ such that 
\begin{align*}
{\bf{P}}(K \,|\, \mathcal{M}_\delta)(\omega)=\E(\chi_{K}(U) \,|\, \mathcal{M}_\delta)(\omega)\quad \text{a.s.}
\end{align*}
Moreover, by applying \cite[Theorem 10.2.1]{Dudley1989} to the joint distribution of $(U,M_\delta)$ we see that there are such functions 
\begin{align*}
(m,K)\mapsto \mathbb{P}_m(K)=:\mathbb{P}\big(\{U\in K\} \,|\, M_\delta=m\big),
\end{align*}
defined for $m\in H^{-s}$ and $K\in \mathcal{B}(H^{\zeta_1})$, that 
\begin{align}\label{conditional}
\mathbb{P}_{M_\delta(\omega)}(K)={\bf{P}}(K \,|\, \mathcal{M}_\delta)(\omega)\quad \text{a.s.}
\end{align}
Using (\ref{credibleset_1}) and (\ref{conditional}) we see that 
\begin{align*}
\mathbb{P}\Big(\big\{U\in B_{\zeta_1}(T_\delta(m),R(\delta))\big\} \,|\, M_\delta=m \Big)
= \mu_\delta\big(B_{\zeta_1}(0,R(\delta)) \big).
\end{align*}
Note that the right hand side is in fact independent of $m$ and depends only on $\delta$.
Next we will give a theorem for the credible sets   
\begin{align*}
\mathbb{P}\Big(\big\{V_\delta\in B_{\zeta_1}(U_\delta,R(\delta))\big\} \Big) 
& = \mathbb{P}\Big(\big\{U\in B_{\zeta_1}(T_\delta(m_\delta),R(\delta))\big\} \,|\, M_\delta=m_\delta \Big).
\end{align*}

\begin{theorem}  
Let $U$, $\Noise$, $M_\delta$ and $A$ be defined as in Theorem \ref{theorem:main} and assume that $A^*:\mathcal{D}'(N)\to\mathcal{D}'(N)$ is invertible.   

Let $\mathcal M_\delta=\sigma(M_\delta)$ and $\mathcal{U}=\sigma(U)$ be the $\sigma$-algebras
generated by the random variable $M_\delta$ and $U$ respectively.
Then the posterior distribution of the random variable $U$ with respect to the $\sigma$-algebra $\mathcal M_\delta$ can be given in terms of function 
\begin{align*}
(m,K)\mapsto \mathbb{P}\big(\{U\in K\} \,|\, M_\delta=m\big),
\end{align*}
where $m\in H^{-s}(N)$ and $K\in \mathcal{B}(H^{\zeta_1})$, cf. (\ref{conditional}).

Take $\zeta_1 < \ntau+t-t_0$ and $\alpha<\gamma/2$. Then if $R(\delta)=C_1\delta^\alpha$ we have the following contraction:   
\ba 
\mathbb{P}\Big(\big\{U\in B_{\zeta_1}(T_\delta(m_\delta),R(\delta))\big\} \,|\, M_\delta=m_\delta \Big)
\geq 1-C\delta^{\gamma-2\alpha} \to 1
\ea
when $\delta\to0$.  
The speed of contraction depends on $\zeta_1$:
\begin{itemize}
\item[(i)] If $\zeta_1\leq -s-t_0$ then $\gamma=2(t+r)/(t_0+r)$.
\item[(ii)] If $-s-t_0\leq \zeta_1<\ntau+t-t_0 $ then $\gamma=2(\ntau+t-t_0-\zeta_1)/(t_0+r)$.
\end{itemize}   
\end{theorem}

\noindent\textit{Proof.}
We use below $R=R(\delta)=C_1\delta^\alpha$ with some $\alpha>0$
and denote
$$
p_\delta=1-\mu_\delta(B_{\zeta_1}(0,R(\delta))).
$$
To study what happens to $p_\delta$ we first notice that   
\begin{align}\label{lowerbound}
\begin{split}
\E\big(\|W_\delta\|_{H^{\zeta_1}}^2\big) &= \int_{H^{\zeta_1}}\|w\|_{H^{\zeta_1}}^2d\mu_\delta(w)\\
&\geq \int_{\|w\|_{H^{\zeta_1}}^2>(R(\delta))^2}\|w\|_{H^{\zeta_1}}^2d\mu_\delta(w)\\
&\geq R(\delta)^2\int_{\|w\|_{H^{\zeta_1}}^2>(R(\delta))^2}d\mu_\delta(w)\\
&= (R(\delta))^2p_\delta.
\end{split}
\end{align}

Next we will prove that 
\ba 
\E\big(\|W_\delta\|_{H^{\zeta_1}}^2\big) = Tr_{H^{\zeta_1}\to H^{\zeta_1}} B_{W_\delta} \leq C\delta^\gamma
\ea
with some $\gamma$ that depends on $\zeta_1$. Above we use the definition $B_{W_\delta}=C_{W_\delta}(I-\Delta)^{\zeta_1}$ where $C_{W_\delta}:H^{-\zeta_1}\to H^{\zeta_1}$ and 
\ba 
\langle C_{W_\delta}\phi,\psi\rangle_{H^{\zeta_1}\times H^{-\zeta_1}} = \E\Big(\langle {W_\delta},\phi \rangle_{H^{\zeta_1}\times H^{-\zeta_1}}\langle {W_\delta},\psi \rangle_{H^{\zeta_1}\times H^{-\zeta_1}}\Big).
\ea
As noted before when $A^*$ is invertible we can write   
\begin{align}\label{def C_h}
C_{W_\delta}=F_\lambda^{-1}(A^*A)^{-1}
\end{align}
where $F_\lambda=( A^*A)^{-1}C_U^{-1}+\lambda I $ and $\lambda=\delta^{-2}$.

We want to estimate
\ba 
\E\big(\|{W_\delta}\|_{H^{\zeta_1}}^2\big)=\E\big(\|(I-\Delta)^{\zeta_1/2}{W_\delta}\|_{L^2}^2\big).
\ea
Let us define 
\ba 
{W_{\delta,\zeta_1}} = (I-\Delta)^{\zeta_1/2}{W_\delta}: \Omega \to L^2(N).
\ea
We can write the covariance operator of ${W_{\delta,\zeta_1}}$
\ba 
C_{W_{\delta,\zeta_1}}=(I-\Delta)^{\zeta_1/2}C_{W_\delta}(I-\Delta)^{\zeta_1/2}:L^2\to L^2.
\ea
Note that in $L^2$ we have $B_{{W_{\delta,\zeta_1}}}=C_{{W_{\delta,\zeta_1}}}$.
Now we get for the trace
\ba 
Tr_{L^2\to L^2}((I-\Delta)^{\zeta_1/2}C_{W_\delta}(I-\Delta)^{\zeta_1/2}) &=& Tr_{L^2\to L^2}((I-\Delta)^{\zeta_1}C_{W_\delta})\\
&=& Tr_{L^2\to L^2}((I-\Delta)^{-s}(I-\Delta)^{\zeta_1+s}C_{W_\delta})\\
&\leq& Tr_{L^2\to L^2}((I-\Delta)^{-s})\|(I-\Delta)^{\zeta_1+s}C_{W_\delta}\|_{L^2\to L^2}.
\ea

Using (\ref{def C_h}) we get
\ba 
\|(I-\Delta)^{\zeta_1+s}C_{W_\delta}\|_{L^2\to L^2} &=& \|(I-\Delta)^{\zeta_1+s}F_\lambda^{-1}(A^*A)^{-1}\|_{0,0}\\
&\leq&  \|(I-\Delta)^{\zeta_1+s}\|_{2(\zeta_1+s),0}\|F_\lambda^{-1}\|_{-2t_0,2(\zeta_1+s)}\|(A^*A)^{-1}\|_{0,-2t_0}\\
&\leq& C\|F_\lambda^{-1}\|_{-2t_0,2(\zeta_1+s)}.
\ea
Above $F_\lambda^{-1}=(( A^*A)^{-1}C_U^{-1}+\lambda)^{-1}\in \Psi_p^m$ where $m=-2(t+r)$ and $p=2(t_0+r)$.
We want to use the norm estimates (\ref{Shubin1}) and (\ref{Shubin2}) so we write
\ba 
\|F_\lambda^{-1}\|_{-2t_0,2(\zeta_1+s)}=\|F_\lambda^{-1}\|_{-2t_0,-2t_0-\ell}
\ea
where $\ell=-2(s+\zeta_1+t_0)$. To use the norm estimates we need to assume $\ell\geq m$, that is, $\zeta_1\leq r-s+t-t_0=\ntau+t-t_0$.

First we assume that $\ell\geq0$ which is true when $\zeta\leq -s-t_0$. Then 
\ba 
\|F_\lambda^{-1}\|_{-2t_0,-2t_0-\ell} &\leq& C(1+\lambda^{\frac{1}{2(t_0+r)}})^{-2(t+r)}
\ea
and $\|F_\lambda^{-1}\|_{-2t_0,-2t_0-\ell}\to 0$ when $\lambda\to\infty$ with all $t,t_0>0$ and $r\geq 0$.

Next we assume $\ell\leq 0$. Then for $-s-t_0\leq \zeta_1\leq \ntau+t-t_0$ we get
\ba 
\|F_\lambda^{-1}\|_{-2t_0,-2t_0-\ell} &\leq& C(1+\lambda^{\frac{1}{2(t_0+r)}})^{2(s+\zeta_1+t_0)-2(t+r)}.
\ea
Now  $\|F_\lambda^{-1}\|_{-2t_0,-2t_0-\ell}\to 0$ when $\lambda\to\infty$ if $\zeta_1< \ntau+t-t_0$. 

We have proven that 
\ba 
\E\big(\|{W_\delta}\|_{H^{\zeta_1}}^2\big)  \leq C\delta^\gamma
\ea
where $\gamma = 2(t+r)/(t_0+r)$ if $\zeta_1\leq -s-t_0$ and $\gamma= 2(\ntau+t-t_o-\zeta_1)/(t_0+r)$ if $-s-t_0\leq \zeta< \ntau+t-t_0$. Hence using  the above estimate and (\ref{lowerbound}) we see that 
\ba 
p_\delta\leq \frac{C\delta^\gamma}{(C_1\delta^\alpha)^2}=C\delta^{\gamma-2\alpha}.
\ea
Above we have to assume $\alpha<\gamma/2$ to have convergence $p_\delta\to 0 $ when $\delta\to 0$.

Finally, since we have denoted $u_\delta=T_\delta(m_\delta)$, we can conclude that with above choices for $R(\delta)$, $\gamma$ and $\alpha$
\begin{align*}
\mathbb{P}\Big(\big\{V_\delta\in B_{\zeta_1}(U_\delta,R(\delta))\big\} \Big) 
& = \mathbb{P}\Big(\big\{U\in B_{\zeta_1}(T_\delta(m_\delta),R(\delta))\big\} \,|\, M_\delta=m_\delta \Big)\\
& \geq 1-C\delta^{\gamma-2\alpha} \to 1 
\end{align*}
when $\delta\to 0$. 
\qed\\

\subsection{Discussion}

Above, we have considered in the  frequentist setting the case when the solution $u^\dagger$ is an element of $H^{\tau}(N)$ with $\tau>0$ and studied in Theorems 2 and 3 the convergence of the MAP estimators and the contraction of the posterior distribution in  $L^2(N)$.  

In the Bayesian setting we have examined the case when the solution is a realisation of the random variable $U$.  In Theorems 1 and 4 we have studied the  the convergence of the MAP estimators and the contraction of the posterior distribution in  $H^\zeta(N)$ with various values of $\zeta$.  

In classical regularisation theory for linear inverse problems, one is usually interested in the convergence of the optimisers of the minimisation problem (2.5) to the true solution in the space $H^r(N)$, $r>\tau+d/2$, as the noise level $\delta$ goes to zero. This gives restrictions to the measurement noise that can be considered. Summarising, our above statistical considerations concern the case where unknown and the noise are significantly less smooth than in the standard setting of the regularisation theory. In the recent regularisation theory inverse problems where the direct map $A$ and the regularisation term are non-linear have been studied extensively.  It is interesting to ask how our analysis on  the contraction of the posterior distribution could be generalised for such non-linear inverse problems that corresponds to  non-Gaussian statistical problems.

\appendix
\section{Some examples of hypoelliptic operators} \label{ap:hypo}

A linear partial differential operator $P$ is hypoelliptic if for every distribution $u$ such that $P(u)$ is $C^\infty$ smooth also $u$ is $C^\infty$. Every elliptic operator with smooths coefficients is hypoelliptic. The heat operator 
\begin{equation*}
Pu(x,t) = \partial_t u-k\Delta_x u,\quad (x,t)\in\R^d\times\R
\end{equation*}
and Kolmogorov operator \cite{kolmogoroff1934,Hormander1967}
\begin{align}\label{Kolmogorov}
Pu(x,y,t) = \partial_{xx}u+x\partial_y u-\partial_t u,\quad (x,y,t)\in\R^3
\end{align}
are examples of operators that are hypoelliptic but not elliptic. General Kolmogorov type hypoelliptic diffusion operators are used e.g. in the theory of kinetic equations, statistical physics and mathematical finance \cite{helffer2005,pascucci2011}.

The fact that (\ref{Kolmogorov}) is hypoelliptic follows from H\"{o}rmander's theorem on hypoelliptic PDEs. Let $(X_0,X_1,\dots,X_p)$ be real $C^\infty$ vector fields in the $d$ dimensional manifold $N$. If $X$ and $Y$ are two vector fields we define the bracket of $X$ and $Y$ by
\begin{align*}
[X,Y]f=X(Yf)-Y(Xf).
\end{align*}
Note that $[X,Y]$ is a new vector field.

\begin{definition}[H\"{o}rmander Condition]\label{hormander}
We say that the H\"{o}rmander condition is satisfied if the real $C^\infty$ vector fields $(X_0,X_1,\dots,X_p)$, $p\leq n$, in the manifold $N$ generate a Lie algebra of rank $n=\dim N$ at every point $x\in N$.
\end{definition}
This means that the vector fields
\begin{align*}
X_j,[X_{j_1},X_{j_2}],[X_{j_1},[X_{j_2},X_{j_3}]],\dots
\end{align*}
span a space that has the same dimension n as the manifold $N$ at every point $x\in N$. Now we can formulate H\"{o}rmander's classical theorem \cite{Hormander1994}.

\begin{theorem}[H\"{o}rmander's theorem]\label{hormander}
The operator 
\begin{align*}
P = \sum_{j=1}^p X_j^2+X_0
\end{align*}
defined on $n$ dimensional manifold $N$
is hypoelliptic if the vector fields $(X_0,X_1,\dots,X_p)$ satisfy the H\"{o}rmander condition.
\end{theorem}

By writing 
\begin{align*}
X_1 = \partial_x, \quad X_0=x\partial_y-\partial_t, \quad [X_1,X_0]=\partial_y
\end{align*}
we see that (\ref{Kolmogorov}) is indeed hypoelliptic in $\R^3$. 
Next we will give another, important example of vector fields satisfying H\"{o}rmander's condition 

\begin{example}
Let us study Heisenberg group $\mathbb{H}$. Let $\Gamma$ be a discrete subgroup of Heisenberg group $\mathbb{H}$ such that $\mathbb{H} / \Gamma$ is compact,see e.g. \cite{folland2004}
The orthonormal frame on $\R^3$ is given by the Lie vector fields
\begin{align*}
X &= \partial_x+\frac{1}{2}y\partial_t\\
Y &= \partial_y-\frac{1}{2}x\partial_t\\
Z &= \partial_t
\end{align*}
We can easily see that 
\begin{align*}
[X,Y] = Z.
\end{align*}
Using Theorem \ref{hormander} we get that the sub-Laplacian 
\begin{align*}
P = \frac{1}{2} (X^2+Y^2)
\end{align*}
on $\mathbb{H}/\Gamma$ is hypoelliptic.
\end{example}

\begin{example}
One example of hypoelliptic inverse problem is the heat equation on a compact manifold $N=M\times\R$ where $M$ is a closed two-dimensional manifold. Note that in this paper we have considered the problem on a compact manifold, that is, our results are applicable in the  case when the equation is periodic in time. We are interested in solving the heat sources $U(x,t)$ from the noisy measurements $M_\delta(x,t)$ of temperature $T(x,t)$, that is we want to solve $U$ from 
\begin{align}\label{heateqex}
(\partial_t-\Delta_x)T(x,t) & = U(x,t), \\
M_\delta(x,t) & = (\partial_t-\Delta_x)^{-1}U(x,t) +\delta\Noise.
\end{align}
The operator $A=(\partial_t-\Delta_x)^{-1}$ is not elliptic but it is hypoelliptic of type (1,2).

Such situations arise in non-invasive monitoring. Consider, for example, using a thermal camera to record video footage of a car with engine running. Let us model the metal surface of the car as a compact and closed two-dimensional manifold $M$. The running engine produces heat which we observe in the video data. The temperature on the car surface is modelled as the solution $T(x,t)$ defined on $N=M\times\R$. Equation \eqref{heateqex} describes the conduction of heat along the car surface. The effect of the engine is simply modelled as the heat source term $U(x,t)$; recovering $U$ will provide information about the state of the engine. 
\end{example}

\section{Computational example}\label{ap: computational}

Since the operator $A$ does not have a continuous inverse operator $L^2\to L^2$,  the condition number of the matrix approximation $\finA$ of the operator $A$ grows when the discretisation is refined.
This is the very reason why regularisation is need in the (numerical) solutions of the inverse problems.

Next we demonstrate the above results numerically and 
 consider two-dimensional deblurring problem {\newtext  on $\T^2$,} 
\begin{equation*}
M = AU+\delta\Noise,
\end{equation*}
where $\Noise(\omega)\in H^{-s}$, $s>1$ a.s. is normalised white noise and $A$ is elliptic operator, smoothing of order $2$,
\[(Au)(x)=\mathcal{F}^{-1}\big((1+|n|^2)^{-1}(\mathcal{F}u)(n)\big)(x).\]
The true solution $u^\dagger\in H^1(\T^2)$, see Subsection  \ref{sec: frequentist} for the frequentist interpretation, is a piecewise linear function presented in Figure \ref{fig: hat}. We choose a priori distribution $N(0,C_U)$ where $C_U=(I-\Delta)^{-1}$. Then a random draw $U(\omega)$ from the prior distribution belongs to $H^{\ntau}$, where $\ntau<0$, with probability one. The Cameron-Martin space $Z_U$ of a measurable mapping $U:\Omega\to X$ is defined by 
\begin{align*}
Z_U=\big\{\phi\in X\, \big|\, \|\phi\|_{Z_U}^2=\langle C_U^{-1}\phi,\phi\rangle_{X^*\times X}<\infty\big\}.
\end{align*}
Cameron-Martin space can also be defined as 
\begin{align*}
Z_U=\bigcap \big\{Y \,\big|\, Y\subset H^{\ntau}\, \text{linear subspace,}\, \mathbb{P}(\{U\in Y\})=1\big\}.
\end{align*}
The approximated solutions $U_\delta$ belongs to $Z_U$ and with the chosen a priori distribution we have $Z_U=H^1(\T^2)$.

\begin{figure}[h]
\centering
\begin{subfigure}[b]{0.5\textwidth}
  \centering
 \includegraphics[height=5.2cm]{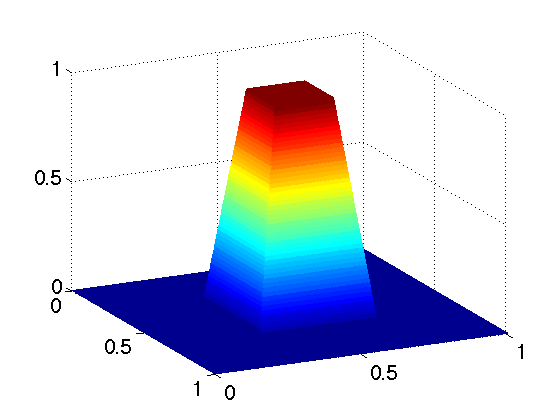}
\end{subfigure}%
\begin{subfigure}[b]{0.5\textwidth}
  \centering
  \includegraphics[height=5.2cm]{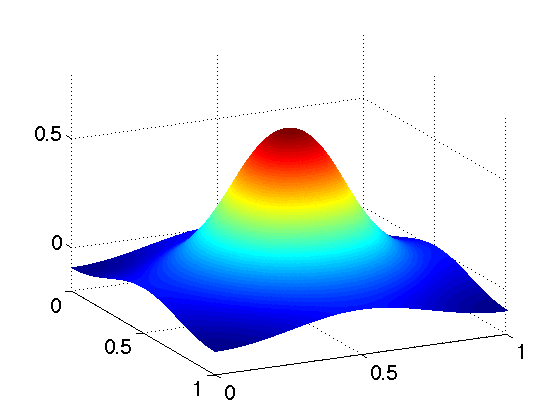}
\end{subfigure}
\caption{On the left the original piecewise linear function $u^\dagger$. On the right side the noiseless data $m^\dagger=Au^\dagger$.}
\label{fig: hat}
\end{figure}

Solving $u$ from $Au(x)=m(x)$ corresponds to the solution of ordinary differential equation $(1-\partial_x^2)m(x)=u(x)$ so $A$ can be thought e.g. as a blurring operator. 

The approximated solution to the problem is
\ba 
u_\delta^\dagger=( A^*A+ \delta^{2}(I-\Delta))^{-1}A^*m^\dagger.
\ea
We get from Subsection \ref{sec: frequentist} that 
\ba 
\lim_{\delta \to 0} \E_u\| u^\dagger-U_\delta^\dagger \|_{H^{\zeta}} = 0
\ea
when $\zeta<\ntau<0$. This behaviour can be seen even in numerical simulations when the discretisation is fine enough, see Figure \ref{fig: norms}. In Figure \ref{fig: comparison} we have compared the expected convergence rates given in formula (\ref{converge_2}) in Theorem \ref{theorem:main} to the computational convergence rates. 
In the numerical simulations in Figures \ref{fig: norms} and \ref{fig: comparison} we see that for the test case presented in Figure \ref{fig: hat} the convergence $u_\delta^\dagger\to u^\dagger$ in different Sobolev spaces follows well the mean convergence predicted by Theorem \ref{theorem:main}. 
\\

\begin{figure}[h!]
\begin{center}
\includegraphics[trim={5cm 0 0 0},clip,height=7cm]{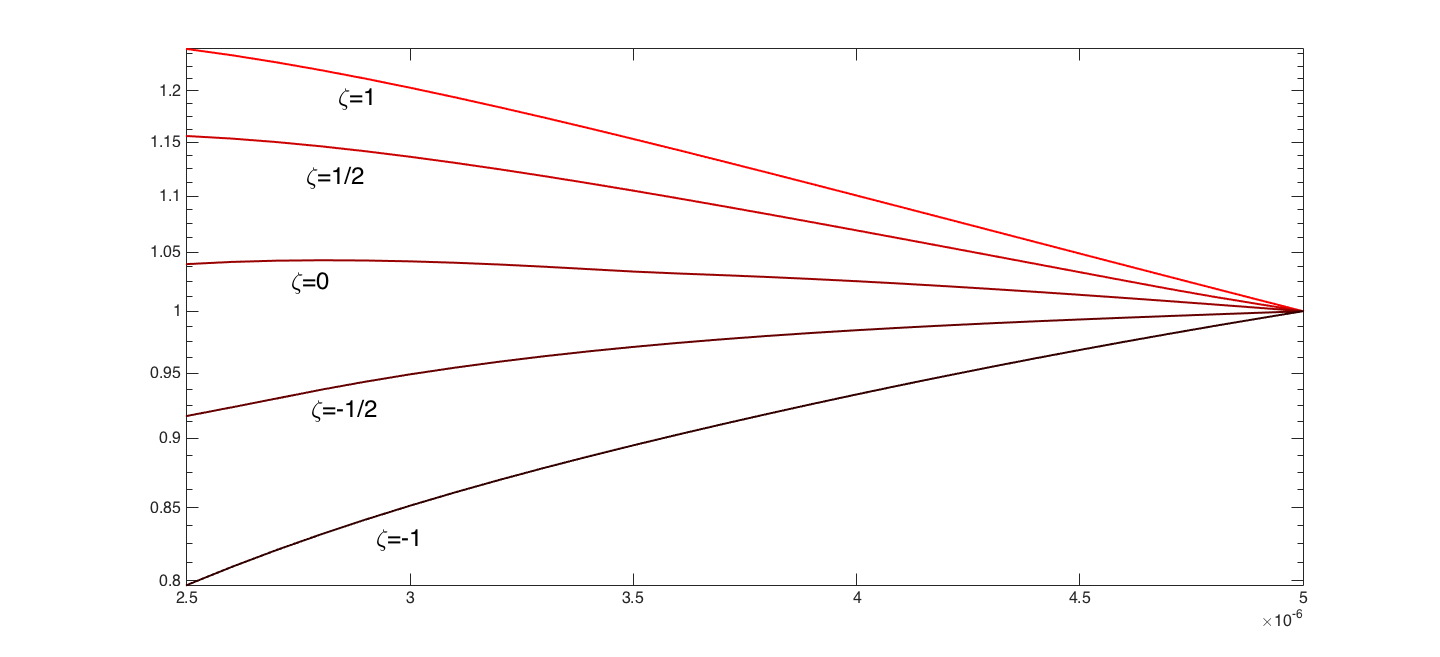}  
\caption{\label{fig: norms}Normalised errors $c(\zeta)\|u^\dagger-u_\delta^\dagger\|_{H^{\zeta}(\T^2)}$ in logarithmic scale with different values of $\zeta$. We use normalisation constants $c(\zeta)=1/\|u^\dagger-u_{5\cdot 10^{-6}}^\dagger\|_{H^{\zeta}}$. We observe that $u_\delta^\dagger\in H^1$ does not converge to $u^\dagger\in H^1$ in $H^{\zeta}$ when $\zeta>0$.}
\end{center}
\end{figure}

\begin{figure}[h!] 
\begin{center}
\includegraphics[trim={5cm 0 0 0},clip,height=7cm]{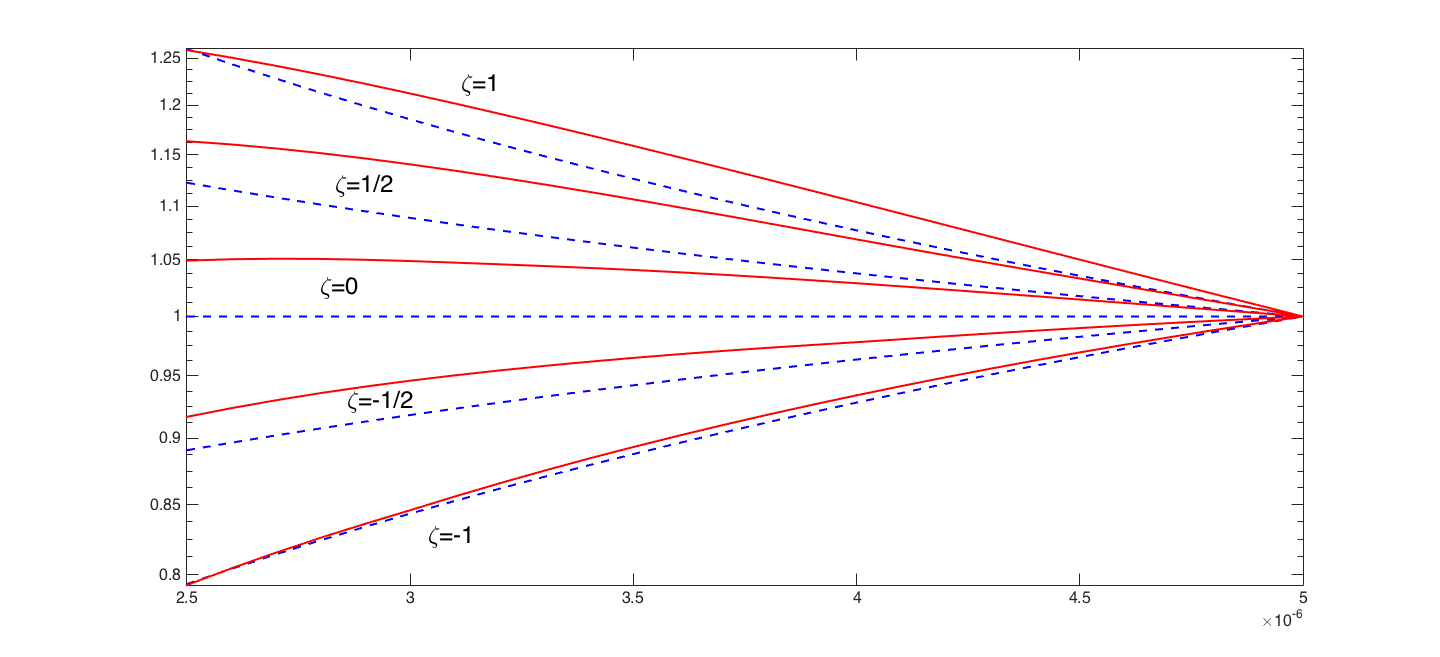}
\caption{\label{fig: comparison}Normalised errors $c(\zeta)\|u^\dagger-u_\delta^\dagger\|_{H^{\zeta}(\T^2)}$ in logarithmic scale with different values of $\zeta$. 
The normalized  bounds (\ref{converge_2}) given in Theorem \ref{theorem:main} for the expectations 
$c_1(\zeta)\E\|U-U_\delta\|_{H^{\zeta}(\T^2)}$
are plotted with dashed lines.
The normalized errors 
$c(\zeta)\|u^\dagger-u_\delta^\dagger\|_{H^{\zeta}(\T^2)}$, for the example $u^\dagger$ given in Figure \ref{fig: hat}, are plotted with solid lines. }
\end{center}
\end{figure}

\medskip
\textbf{Acknowledgements.}
We would like to thank Petteri Piiroinen for helpful discussions. 
This work was supported by the Finnish Centre of Excellence in Inverse Problems Research 2012-2017 (Academy of Finland CoE-project
284715). In addition, H.K. was supported by Emil Aaltonen Foundation and EQUIP, grant EP/K034154/1, M.L. was supported by Academy of Finland, grant
273979, and S.S. was supported by Academy of Finland, project 141094.

\newpage
\bibliographystyle{siam}
\bibliography{Inverse_problems_references}

\end{document}